\documentclass[11pt, twoside]{amsart}
\usepackage[dvipsnames]{xcolor}
\usepackage{hyperref}
\hypersetup{colorlinks=true, linkcolor=Blue, citecolor=Maroon}
\usepackage{amsthm, amsmath, amscd, amssymb,centernot,txfonts}
\usepackage{tikz-cd}
\usepackage{float}
\usepackage{multirow}
\usepackage{lipsum}
\usepackage{multicol}
\usepackage{amsmath}
\usepackage{amssymb}

\usepackage[normalem]{ulem}
\usepackage{amsfonts}
\usepackage{mathrsfs}
\usepackage{parskip}
\usepackage[all]{xy}
\usepackage[left=2.5cm,top=2.5cm,bottom=3cm,right=2.5cm]{geometry}
\usepackage{dirtytalk}
\usepackage{verbatim}
\usepackage{mathtools}
\usetikzlibrary{shapes,arrows}
\usepackage{etoolbox}
\usepackage{dynkin-diagrams}

\usepackage{enumitem}
\usepackage{chngpage}

\usepackage{adjustbox}
\usepackage[utf8]{inputenc}
\usepackage{fourier} 
\usepackage{array}
\usepackage{makecell}

\usepackage[linktocpage=true]{}

\usepackage{color, colortbl}
\definecolor{LightCyan}{rgb}{0.88,1,1}

\usepackage[normalem]{ulem}

\usepackage[first=0,last=9]{lcg}

\DeclareMathOperator*\lowlim{\underline{lim}}
\DeclareMathOperator*\uplim{\overline{lim}}

\newcommand{\stkout}[1]{\ifmmode\text{\sout{\ensuremath{#1}}}\else\sout{#1}\fi}

\theoremstyle{plain}
\numberwithin{equation}{section}
\newtheorem{theorem}{Theorem}[section]
\newtheorem{proposition}[theorem]{Proposition}
\newtheorem{lemma}[theorem]{Lemma}
\newtheorem{corollary}[theorem]{Corollary}

\newtheorem{question}[theorem]{Question}

\newtheorem{set-up}[theorem]{Set-up}
\newtheorem{observation}[theorem]{Observation}

\theoremstyle{definition}
\newtheorem{remark}[theorem]{Remark}

\newtheorem{example}[theorem]{Example}
\newtheorem{definition}[theorem]{Definition}

\newcommand*{\QEDB}{\hfill\ensuremath{\square}}

\tikzstyle{decision} = [diamond, draw, , 
    text width=4.5em, text badly centered, node distance=3cm, inner sep=0pt]
\tikzstyle{block} = [rectangle, draw, , 
    text width=10em, text centered, rounded corners, minimum height=2em]
\tikzstyle{block1} = [rectangle, draw, , 
    text width=5em, text centered, rounded corners, minimum height=2em]    
\tikzstyle{line} = [draw, -latex']
\tikzstyle{cloud} = [draw, ellipse,, node distance=3cm,
    minimum height=2em]

\begin{document}
\title[Continuous CM--regularity and generic vanishing]{Continuous CM--regularity and generic vanishing}

\author[]{}
\address{Department of Mathematics, Faculty of Natural Science and Technology,
Okayama University,
Okayama, Japan}
\email{ito-atsushi@okayama-u.ac.jp}

\author[]{Debaditya Raychaudhury, \hspace{3pt} with an appendix by Atsushi Ito}
\address{Department of Mathematics, University of Toronto, Toronto, Canada}
\email{debaditya.raychaudhury@utoronto.ca}
\address{Current address of D. Raychaudhury: Department of Mathematics, University of Arizona,
Tucson, USA}
\email{draychaudhury@math.arizona.edu, rcdeba@gmail.com}


\maketitle

\vspace{-10pt}

\begin{abstract}
We study the continuous CM--regularity of torsion--free coherent sheaves on polarized irregular smooth projective varieties $(X,\mathcal{O}_X(1))$, and its relation with the theory of generic vanishing. This continuous variant of the Castelnuovo--Mumford regularity was introduced by Mustopa, and he raised the question whether a continuously $1$--regular such sheaf $\mathcal{F}$ is GV. Here we answer the question in the affirmative for many pairs $(X,\mathcal{O}_X(1))$ which includes the case of any polarized abelian variety. Moreover, for these pairs, we show that if $\mathcal{F}$ is continuously $k$--regular for some integer $1\leq k\leq \dim X$, then $\mathcal{F}$ is a GV$_{-(k-1)}$ sheaf. Further, we extend the notion of continuous CM--regularity to a real valued function on the $\mathbb{Q}$--twisted bundles on polarized abelian varieties $(X,\mathcal{O}_X(1))$, and we show that this function can be extended to a continuous function on $N^1(X)_{\mathbb{R}}$. We also provide syzygetic consequences of our results for $\mathcal{O}_{\mathbb{P}(\mathcal{E})}(1)$ on $\mathbb{P}(\mathcal{E})$ associated to a $0$--regular bundle $\mathcal{E}$ on polarized abelian varieties. In particular, we show that $\mathcal{O}_{\mathbb{P}(\mathcal{E})}(1)$ satisfies $N_p$ property if the base--point freeness threshold of the class of $\mathcal{O}_X(1)$ in $N^1(X)$ is less than $\frac{1}{p+2}$. This result is obtained using a theorem in the \hyperref[appendix]{Appendix A} written by Atsushi Ito. 
\end{abstract}

\section{Introduction} 

Given a smooth projective variety $X\subseteq\mathbb{P}^n$, it is well--known that the geometry of the embedding is reflected by the coherent sheaves on $X$ with suitable positivity properties. One of the most fundamental notion of such properties that governs the complexity of a sheaf $\mathcal{F}$ is given by its Castelnuovo--Mumford (CM) regularity with respect to the pair $(X,\mathcal{O}_X(1))$. 

\smallskip

The purpose of this article is to study a variant of CM--regularity that was introduced by Mustopa in \cite{Mus} for polarized irregular varieties $(X,\mathcal{O}_X(1))$. It is defined as follows: $\mathcal{F}$ is {\it continuously $k$--regular} if the cohomological support loci $V^i(\mathcal{F}(k-i))\neq\emptyset$ for $i\geq 1$ (strictly speaking, this definition is slightly more general from that in \cite{Mus} as we are not assuming global generation of $\mathcal{O}_X(1)$). The structures of these cohomological support loci are of great importance in the study of the geometry of irregular varieties. An important notion in this topic is the notion of generic vanishing (GV for short) introduced by Green and Lazarsfeld in the pioneering works \cite{GL1}, \cite{GL2}. Fundamental contributions from Hacon (\cite{Hac}), and Pareschi--Popa (\cite{PaPo1} -- \cite{PaPo4}) through derived category approach and Fourier--Mukai functors led to subsequent developments of the theory of generic vanishing. 

\smallskip

Turning to details, a coherent sheaf $\mathcal{F}$ on $X$ is said to be {\it GV} if codim$V^i(\mathcal{F})\geq i$ for all $i$. This property is intimately related to the positivity of $\mathcal{F}$, in particular, a GV sheaf on an abelian variety is nef. In \cite{Mus}, Mustopa asked the following question.

\begin{question}\label{mq}
(\cite[Question $(*)$]{Mus}) Let $X$ be a smooth projective variety of dimension $d\geq 1$ and let $\mathcal{O}_X(1)$ be an ample and globally generated line bundle on $X$. Let $\mathcal{F}$ be a torsion--free coherent sheaf on $X$. If $\mathcal{F}$ is continuously $1$--regular for $(X,\mathcal{O}_X(1))$, is $\mathcal{F}$ a GV sheaf?
\end{question}

The question above was motivated by Beauville's construction in \cite{Beauab} of rank 2 Ulrich bundles on abelian surfaces $(X,\mathcal{O}_X(1))$. It turns out that for these bundles $\mathcal{E}$, $\mathcal{E}(-1)$ is indeed GV. It is easy to see that the answer of Question \ref{mq} is affirmative for polarized curves. It was shown in \cite{Mus} that the answer of the question is also affirmative for:
\begin{itemize}[leftmargin=0.15in]
    \item[--] a large class of polarized surfaces that includes the case of any polarized abelian surface ({\it loc. cit}. Theorem B, Corollary C),
    \item[--] certain polarizations on Cartesian and symmetric products of curves ({\it loc. cit.} Propositions 3.1, 3.2),
    \item[--] some scrollar embeddings of ruled threefolds over a curve ({\it loc. cit.} Proposition 4.1).
\end{itemize}

\smallskip

Continuous CM--regularity for semihomogeneous bundles on abelian varieties has been studied by K\"uronya and Mustopa in \cite{KM}, and later by Grieve in \cite{Gre21}. In particular, K\"uronya and Mustopa showed in \cite{KM}, that if $\mathcal{E}$ is a semihomogeneous bundle on abelian variety $(X,\mathcal{O}_X(1))$ of dimension $g$ as in the set--up of Question \ref{mq}, and moreover if $c_1(\mathcal{E})$ is a rational multiple of $c_1(\mathcal{O}_X(1))$, then even more is true, namely $\mathcal{E}(1-g)$ is GV. In \cite{Gre21}, Grieve established a description of continuous CM--regularity of semihomogeneous bundles on abelian varieties. This description was in terms of a normalized polynomial function studied in \cite{Gri17}, and obtained via the Wedderburn decomposition of the endomorphism algebra of the abelian variety. A further study of index and generic vanishing theory of simple semihomogeneous bundles was also carried out in \cite{Gre21} building on and refining \cite{Gri14a}, \cite{Gri14b}.

\smallskip

Beside the notion of GV sheaves, a related notion in the generic vanishing theory is that of M--regularity of coherent sheaves. A coherent sheaf $\mathcal{F}$ on $X$ is said to be {\it M--regular} if codim$V^i(\mathcal{F})> i$ for all $i>0$. In this direction, Mustopa asked whether a continuous $0$--regular torsion--free coherent sheaf on a polarized smooth projective variety $(X,\mathcal{O}_X(1))$ with $\mathcal{O}_X(1)$ globally generated is M--regular (\cite[Remark 1.6]{Mus}).

\smallskip

Following this train of thought, it is natural to propose the more general question: given a continuously $k$--regular torsion--free coherent sheaf on a smooth polarized variety $(X,\mathcal{O}_X(1))$ with $\mathcal{O}_X(1)$ globally generated and $1\leq k\leq\dim X$, is it true that codim$V^i(\mathcal{F})\geq i-k+1$ for all $i$? Here we remark that the notion of generic vanishing was generalized in \cite{PaPo4} where Pareschi--Popa called a sheaf $\mathcal{F}$ {\it GV}$_{-k}$ for an integer $k\geq 0$, if codim$(V^i(\mathcal{F}))\geq i-k$ for all $i$. In view of this, here we will devote ourselves to find the answer of the question below which is more general than \cite[Question $(*)$ and Remark 1.6]{Mus}.

\begin{question}\label{newq}
Let $X$ be a smooth projective variety of dimension $d\geq 1$ and let $\mathcal{O}_X(1)$ be an ample and globally generated line bundle on $X$. Let $\mathcal{F}$ be a torsion--free coherent sheaf on $X$. Assume $\mathcal{F}$ is continuously $k$--regular for $(X,\mathcal{O}_X(1))$ for some integer $0\leq k\leq d$.
\begin{itemize}
    \item[(1)] If $1\leq k\leq d$, is $\mathcal{F}$ a GV$_{-(k-1)}$ sheaf?
    \item[(2)] If $k=0$, is $\mathcal{F}$ an M--regular sheaf?
\end{itemize}
\end{question}

The following is the main result of this article that answers the above question in the affirmative for many pairs $(X,\mathcal{O}_X(1))$. We also note that the result below does not require the hypothesis that the polarization $\mathcal{O}_X(1)$ is globally generated.

\vspace{5pt}

\noindent {\bf Theorem A.}\phantomsection\label{newgv}
\textit{Let $(X,H)$ be a polarized smooth projective variety. Assume that there exists a globally generated line bundle $H_1$ on $X$ and an ample line bundle $H_2$ on $\textrm{Alb}(X)$ such that $H=H_1+\textrm{alb}_X^*H_2$, where $\textrm{alb}_X:X\to\textrm{Alb}(X)$ is the Albanese map. Let $\mathcal{F}$ be a torsion--free coherent sheaf on $X$ that is continuously $k$--regular for $(X,H)$ for some integer $0\leq k\leq \dim X$. Then the following statements hold.
\begin{itemize}
    \item[(1)] $V^i(\mathcal{F})=\emptyset$ for $i\geq k+1$.
    \item[(2)] If $k\neq 0$, then $\textrm{codim}(V^k(\mathcal{F}))\geq 1$.
\end{itemize} 
In particular, the answer of Question \ref{newq} is affirmative for the pair $(X,H)$.} 

\vspace{5pt}

Note that the above result answers Question \ref{newq} in the affirmative for any polarized abelian variety. Moreover, it also shows that the question has an affirmative answer in many cases that include for example (a) $X=Y\times A$ where $Y$ is a regular smooth projective variety, $A$ is an abelian variety, and $\mathcal{O}_X(1)$ is ample and globally generated; and (b) $X$ is a projective bundle $\mathbb{P}(\mathcal{E})$ on an abelian variety $A$ associated to an ample and globally generated vector bundle $\mathcal{E}$, and $\mathcal{O}_X(1)=T+F$ where $T$ is the tautological bundle, and $F$ is the pull--back of any ample line bundle on $A$. Observe also that the above theorem proves a stronger statement than what is asked for in Question \ref{newq} (2) when $(X,H)$ satisfies the hypotheses of Theorem \hyperref[newgv]{A}: it shows that continuously $0$-regular sheaves are in fact IT$_0$ (see Definition \ref{defit} and Corollary \ref{gvnef}). We will give other variants of Theorem \hyperref[newgv]{A} in \S \ref{secvar}.

\smallskip

The proof of Theorem \hyperref[newgv]{A} is homological in nature and relies on an inductive argument, but the crucial ingredient of the proof is a relative set--up of continuous CM--regularity that we develop in this article. 

\smallskip

Inspired by the recent development of the {\it cohomological rank function} by Jiang--Pareschi (see \cite{JP}), we further extend the notion of continuous CM--regularity to define a real--valued regularity  $\mathbb{Q}\textrm{-reg}_{\underline{l}}(\mathcal{E}\langle \underline{n}\rangle)$ for $\mathbb{Q}$--twisted bundles $\mathcal{E}\langle \underline{n}\rangle$ where $\underline{l}$ is an ample class in $N^1(X)$ on an abelian variety $X$, and $\underline{n}\in N^1(X)_{\mathbb{Q}}$. For them, we prove the following

\vspace{5pt}

\noindent{\bf Theorem B.}\phantomsection\label{qregthm}
{\it Let $X$ be an abelian variety of dimension $g$, and let $\mathcal{E}$ be a vector bundle on $X$. If $\underline{l}\in N^1(X)$ is a polarization, then the function $\mathbb{Q}\textrm{-reg}_{\underline{l}}(\mathcal{E}\langle - \rangle): N^1(X)_{\mathbb{Q}}\to\mathbb{R}$ sending $N^1(X)_{\mathbb{Q}}\ni\underline{n}\mapsto\mathbb{Q}\textrm{-reg}_{\underline{l}}(\mathcal{E}\langle\underline{n}\rangle)$ can be extended to a continuous function $\mathbb{R}\textrm{-reg}_{\underline{l}}(\mathcal{E}\langle - \rangle):N^1(X)_{\mathbb{R}}\to\mathbb{R}$.}

\vspace{5pt}

We now mention the immediate consequences of Theorem \hyperref[newgv]{A} for continuously $k$--regular torsion--free coherent sheaves $\mathcal{F}$ on polarized abelian varieties $(X,\mathcal{O}_X(1))$. It follows immediately that
\begin{itemize}
    \item (Corollary \ref{v0}) if $k=1$, then $\chi(\mathcal{F})\geq 0$ with equality if and only if $V^0(\mathcal{F})$ is a divisor,
    \item (Corollary \ref{gvnef}) if $k=1$ then $\mathcal{F}$ is nef, and if $k=0$ then it is ample.
\end{itemize}
Let us denote the (usual) CM--regularity by $\textrm{reg}_{\mathcal{O}_X(1)}(-)$. In Corollary \ref{regin} we establish the sub--additivity of CM--regularity for polarized abelian varieties. To be more precise, for any torsion--free coherent sheaves $\mathcal{E}$ and $\mathcal{F}$ on a polarized abelian variety $(X,\mathcal{O}_X(1))$  such that at least one of $\mathcal{E}$ and $\mathcal{F}$ is locally free, we show that the following inequality holds 
    \begin{equation}\label{rege1e2}
        \textrm{reg}_{\mathcal{O}_X(1)}(\mathcal{E}\otimes\mathcal{F})\leq \textrm{reg}_{\mathcal{O}_X(1)}(\mathcal{E})+\textrm{reg}_{\mathcal{O}_X(1)}(\mathcal{F}).
    \end{equation} 
We remark that it was shown earlier by Totaro in \cite[Theorem 3.4]{Tot} that \eqref{rege1e2} holds for arbitrary polarized smooth projective variety when $\mathcal{O}_X(1)$ is very ample and satisfies a certain Koszul hypothesis, see also \cite[Theorem 1.1]{ES} for the case $X=\mathbb{P}^d$.

\smallskip

We also deduce syzygetic consequences from Theorem \hyperref[newgv]{A} for projective bundles on abelian varieties. Before stating our result, we spend a few words on linear series on abelian varieties to set the context. 

\smallskip

It was a conjecture of Lazarsfeld that on a polarized abelian variety $(X,H)$, $tH$ satisfies $N_p$ property if $t\geq p+3$. Lazarsfeld's conjecture was proven by Pareschi in \cite{Par00}. The result was further extended by Pareschi--Popa in \cite{PaPo2} where it was shown that $tH$ satisfies $N_p$ property if $t\geq p+2$ and $|H|$ contains no base--divisor. Further, let us denote the ideal sheaf at the origin of $X$ as $\mathcal{I}_0$, and define 
$$r(H):=\textrm{inf}\left\{c\in\mathbb{Q}\mid \textrm{there exists an effective $\mathbb{Q}$--divisor $D\equiv cH$ such that $\mathcal{J}(X,D)=\mathcal{I}_0$}\right\}$$ 
where $\mathcal{J}(X,D)$ is the multiplier ideal. Lazarsfeld--Pareschi--Popa proved in \cite{LPP} that if $r(H)<\frac{1}{p+2}$ then $H$ satisfies $N_p$ property. Very recently, Jiang--Pareschi defined in \cite{JP} the {\it base--point freeness threshold} $0<\beta(\underline{h})\leq 1$ for a polarization $\underline{h}\in N^1(X)$ on an abelian variety $X$ that can be characterized as follows: $\beta(\underline{h})<x\iff \mathcal{I}_0\langle x\underline{h}\rangle$ is IT$_0$ for $x\in\mathbb{Q}$. It was shown by Caucci in \cite{Cau} that $\beta(\underline{h})\leq r(H)$, and further the above results are unified to shown that $H$ satisfies $N_p$ property if $\beta(\underline{h})<\frac{1}{p+2}$. The study of a more general property $N_p^r$ via base--point freeness threshold was carried out by Ito in \cite{Ito22}. Sharp results on projective normality and higher syzygies of general polarized abelian varieties were also established by Ito in \cite{Ito21} and in \cite{Ito22'}.

\smallskip

The following is our result on syzygies of projective bundles associated to continuously $0$--regular vector bundles on polarized abelian varieties $(X,H)$, which is an immediate consequence of Theorem \hyperref[newgv]{A} and Theorem \ref{thm_appendix} in the \hyperref[appendix]{Appendix A} which is due to Atsushi Ito.

\vspace{5pt}

\phantomsection\label{npmain}
\noindent{\bf Corollary C.} {\it Let $(X,H)$ be a polarized abelian variety. Also let $\mathcal{E}$ be a continuously $0$--regular vector bundle for $(X,H)$. Then $\mathcal{O}_{\mathbb{P}(\mathcal{E})}(1)$ satisfies $N_{p}$ property if $\beta(\underline{h})<\frac{1}{p+2}$ where $\underline{h}$ is the class of $H$ in $N^1(X)$.}

\subsection{Organization} In \S \ref{4}, we recall the preliminaries of (continuous) CM-regularity. In \S \ref{4.5} we first discuss the theory of generic vanishing, and then we proceed to prove Theorem \hyperref[newgv]{A}. At the end of this section, we prove a few variants of Theorem \hyperref[newgv]{A}. We define the $\mathbb{Q}$ CM--regularity \S \ref{4.7} and in this section we prove Theorem \hyperref[qregthm]{B}. Finally \S \ref{seccor} is devoted to the proofs of Corollaries \ref{v0}, \ref{gvnef}, \ref{regin} and Corollary \hyperref[npmain]{C}.  

\subsection{Conventions} We work over the field of complex numbers $\mathbb{C}$. We will tacitly assume that the varieties are irregular and the morphisms are non--constant. We use the additive and multiplicative notation interchangeably for tensor products of line bundles, and the sign ``$\equiv$'' will be used for numerical equivalence. The rest of the notation is standard in algebraic geometry. 

\subsection{Acknowledgements} It is my great pleasure to acknowledge my gratitude to Professor Purnaprajna Bangere, and Professor Angelo Felice Lopez for their constant encouragement and support. In particular, I am extremely grateful to Professor Angelo Felice Lopez for suggesting me to think about projective normality of Ulrich bundles, which is how this project got started. I am indebted to Professor Atsushi Ito for generously sharing his ideas on syzygies of projective bundles on abelian varieties, and for kindly agreeing to write the \hyperref[appendix]{Appendix A} of this article. I sincerely thank Professors Nathan Grieve, V. Kumar Murty, Robert Lazarsfeld, Yusuf Mustopa, and Mihnea Popa for their valuable comments and suggestions on an earlier draft. I also benefited from conversations with Jayan Mukherjee for which I thank him. I am grateful to the anonymous referee for a careful reading of the manuscript and for suggesting several improvements. During the preparation of this work, I was supported by a Simons Postdoctoral Fellowship from the Fields Institute for Research in Mathematical Sciences.

\section{Continuous CM--regularity of coherent sheaves}\label{4}

\subsection{Definition and first properties} This subsection is devoted to the definitions and basic properties of (continuous) CM--regularity. We start with the definition of a partial variant of the usual CM--regularity (cf. \cite[ Definition 1.8.4]{Laz}, and \cite[Lemma 3.2]{Tot}).

\begin{definition}\label{defreg}
{\bf (CM--regularity)} Let $X$ be a smooth projective variety and let $H$ be a line bundle on $X$. Also, let $q,k$ be integers with $q\geq 0$. A coherent sheaf $\mathcal{F}$ on $X$ is called $C_{q,k}$ for $(X,H)$ if $H^{q+i}(\mathcal{F}((k-i)H))=0$ for all integers $i\geq 1$. When $H$ is ample, we will say that $\mathcal{F}$ is {\it $k$--regular for $(X,H)$} for $k\in\mathbb{Z}$ if it is $C_{0,k}$ for $(X,H)$.
\end{definition}

The following important result in the study of partial regularity was proven by Totaro in \cite{Tot}, and is well--known in the set--up of the usual (i.e., non--partial) Castelnuovo--Mumford regularity. 

\begin{lemma}\label{partial}
(\cite[Lemma 3.2]{Tot}) Let $(X,H)$ be a smooth projective variety with $H$ globally generated. Let $\mathcal{F}$ be a coherent sheaf on $X$, and $q,k$ be integers with $q\geq 0$. If $\mathcal{F}$ is $C_{q,k}$ for $(X,H)$ then it is $C_{q,k+m}$ for $(X,H)$ for any integer $m\geq 0$.
\end{lemma}

When $H$ is ample and globally generated, define the {\it regularity} of $\mathcal{F}$ as $$\textrm{reg}_{H}(\mathcal{F}):=\textrm{min}\left\{m\in\mathbb{Z}\mid\textrm{ $\mathcal{F}$ is $m$--regular for $(X,H)$}\right\}.$$

\smallskip

We now recall the definition of the cohomological support loci (cf. \cite[Definition 1.2]{Par'}) that are of fundamental importance in the study of irregular varieties.

\begin{definition} {\bf(Cohomological support loci)} Let $X$ be a smooth projective variety and let $a:X\to A$ be a morphism to an abelian variety $A$. Let $\mathcal{F}$ be a coherent sheaf on $X$. The {\it $i$--th cohomological support locus $V^i_a(\mathcal{F})$ with respect to $a$} for $i\in\mathbb{N}$ is defined as
\begin{equation*}
    V^i_a(\mathcal{F}):=\left\{\zeta\in\textrm{Pic}^0(A)\,|\,h^i(\mathcal{F}\otimes a^*{\zeta})\neq 0\right\}.
\end{equation*}
We will simply write $V^i(\mathcal{F})$ for $V^i_{\textrm{alb}_X}(\mathcal{F})$ where $\textrm{alb}_X:X\to\textrm{Alb}(X)$ is the Albanese morphism of $X$.
\end{definition}
    
As we discussed earlier, continuous CM--regularity is a slightly coarser measure of positivity than CM--regularity, and was introduced by Mustopa in \cite{Mus}. For our purpose, we need to generalize the definition of continuous CM--regularity (cf. \cite[Definition 1.1]{Mus}) to a  relative and partial set--up that we describe next.

\begin{definition}\label{defcontcm}
{\bf(Continuous CM--regularity)} Let $X$ be a smooth projective variety and let $H$ be a line bundle on $X$. Further, let $a:X\to A$ be a morphism to an abelian variety $A$. Let $\mathcal{F}$ be a coherent sheaf on $X$ and let $q,k$ be integers with $q\geq 0$. The sheaf $\mathcal{F}$ is called $C'^a_{q,k}$ for $(X,H)$ if $V_a^{q+i}(\mathcal{F}((k-i)H))\neq \textrm{Pic}^0(A)$ for all integers $i\geq 1$. When $H$ is ample, we will say $\mathcal{F}$ is {\it continuously $k$--regular for $(X,H)$} if it is $C'^{\textrm{alb}_X}_{0,k}$ for $(X,H)$.

\end{definition}

When $H$ is ample and globally generated, we define the {\it continuous regularity} of $\mathcal{F}$ as
$$\textrm{reg}_H^{\textrm{cont}}(\mathcal{F}):=\textrm{min}\left\{m\in\mathbb{Z}\,|\, \forall i>0, V^i(\mathcal{F}((m-i)H))\neq\textrm{Pic}^0(\textrm{Alb}(X))\cong\textrm{Pic}^0(X)\right\}.$$
In general, we have the inequality $\textrm{reg}_H^{\textrm{cont}}(\mathcal{F})\leq\textrm{reg}_H(\mathcal{F})$. However, strict inequalities are possible:

\begin{example}
    Let $H$ be an ample and globally generated line bundle on an abelian variety $X$ of dimension $g$. Then $g=\textrm{reg}_H^{\textrm{cont}}(\mathcal{O}_X)< \textrm{reg}_H(\mathcal{O}_X)=g+1$.
\end{example}

We will use the following fact without any further reference.

\begin{remark}
    Let $a:X\to A$ be a morphism from a smooth projective variety $X$ to an abelian variety $A$. Let $\mathcal{F}$ be a coherent sheaf and $H$ be a line bundle on $X$. It is a consequence of Definition \ref{defreg} (resp. \ref{defcontcm}) that if $\mathcal{F}$ is $C_{q,k}$ (resp. $C'^a_{q,k}$) for $(X,H)$ for integers $q,k\geq 0$, then it is also $C_{q+k,0}$ (resp.  $C'^a_{q+k,0}$). 
\end{remark}

We will see that in practice, it is often useful to work with relative continuous CM--regularity rather than CM--regularity. The following observation (where we use semicontinuity to see that the first three equivalent condition implies the fourth) highlights this and shows that the former property is stable under perturbations by elements of $a^*\textrm{Pic}^0(A)$.

\begin{observation}\label{obs}
(see also \cite[Lemma 1.2]{Mus}) Let $X$ be a smooth projective variety and let $H$ be a line bundle on $X$. Let $a:X\to A$ be a morphism to an abelian variety $A$, and let $\mathcal{F}$ be a coherent sheaf on $X$. Also, let $q,k$ be integers with $q\geq 0$. The following conditions are equivalent:
\begin{itemize}
    \item[(1)] $\mathcal{F}$ is $C'^a_{q,k}$ for $(X,H)$,
    \item[(2)] $\mathcal{F}\otimes a^*\zeta$ is $C'^a_{q,k}$ for $(X,H)$ for some (equivalently for all) $\zeta\in\textrm{Pic}^0(A)$,
    \item[(3)] $\mathcal{F}$ is $C'^a_{q,k}$ for $(X,H+a^*\zeta)$ for some (equivalently for all) $\zeta\in\textrm{Pic}^0(A)$,
    \item[(4)] $\mathcal{F}\otimes a^*\zeta$ is $C_{q,k}$ for $(X,H)$ for some (equivalently for general) $\zeta\in\textrm{Pic}^0(A)$.
\end{itemize}
\end{observation}

Thus, the (partial) continuous CM--regularity is determined by the class of the line bundle $H$ in the N\'eron--Severi group $\textrm{Pic}(X)/\textrm{Pic}^0(X)$. 

\subsection{Behavior of continuous CM--regularity} In this subsection, we study the behavior of continuous CM--regularity of torsion--free coherent sheaves with respect to
restriction and pull--back.

\subsubsection{} We introduce the property $(P_a)$ that is crucial for us since it allows us to produce smooth sections in the appropriate linear series by Bertini's theorem. 

\begin{definition}\label{propstar}
{\bf (Property $(P_a)$)} We will say that a polarized smooth projective variety $(X,H)$ satisfies $(P_a)$ where $a:X\to A$ is a morphism to an abelian variety if for all $\zeta\in\textrm{Pic}^0(A)$, $H+a^*\zeta$ is globally generated.
\end{definition}

The above property is desirable for various geometric reasons aside from the one we stated before, we point out another such reason here. A useful notion for sheaves on irregular varieties is their {\it continuous global generation} property that we will define in Definition \ref{cggprop}. It follows from \cite[Proposition 3.1]{Deb} that if $H$ satisfies Property $(P_{\textrm{alb}_X})$, then $H$ is continuously globally generated.

\smallskip

We now make the following

\begin{observation}\label{obs2}
Let $(X,H)$ be a polarized smooth projective variety and let $a:X\to A$ be a morphism to an abelian variety. Assume that $(X,H)$ satisfies $(P_a)$. Then for any $\zeta\in\textrm{Pic}^0(A)$ and any smooth and irreducible member $Y\in|H+a^*\zeta|$, the pair $(Y,H|_Y)$ satisfies $(P_{a|_Y})$ where $a|_Y:Y\to A$ is the restriction of $a$.
\end{observation}

It is important for us to note that continuous CM--regularity satisfies better restriction properties than ordinary CM--regularity. We highlight this by means of an example which shows that in general, if $\mathcal{E}$ $k$--regular for $(X,H)$, $\mathcal{E}|_Y$ need not be $k$--regular for $(Y,H|_Y)$ if $Y\in|H+\xi|$ where $\xi$ is non--trivial in Pic$^0(X)$. 

\begin{example}
    Let $X$ be an abelian surface, and let $H$ be an ample and globally generated line bundle. Fix a line bundle $0\neq\xi\in \textrm{Pic}^0(X)$ and observe that $H+\xi$ is globally generated (this fact has been pointed out to me by the referee, whom I thank). Indeed, consider the isogeny $\varphi_H: X\to \hat{X}=\textrm{Pic}^0(X)$ that sends $x\in X$ to $t_x^*H\otimes H^{\otimes -1}$ where $t_x:X\to X$ is translation by $x$. Since $\varphi_H$ is an isogeny, in particular surjective, there exists $x\in X$ such that $H+\xi=t_x^*H$ whence the global generation of $H+\xi$ follows. Thus, $H+\xi$ is ample and globally generated, and consequently there exists a smooth curve $Y\in |H+\xi|$. It is evident that $2H+\xi$ is $0$--regular for $(X,H)$. However, we claim that $(2H+\xi)|_Y$ is not $0$--regular for $(Y,H|_Y)$. To see this, observe that $(H+\xi)|_Y=K_Y$ by adjunction, whence $h^1((H+\xi)|_Y)\neq 0$ by Serre duality. 
\end{example}

However, for the continuous variant of CM--regularity, we have 

\begin{lemma}\label{contres}
Let $(X,H)$ be a polarized smooth projective variety and let $a:X\to A$ be a morphism to an abelian variety $A$. Assume that $(X,H)$ satisfies $(P_a)$. Let $\mathcal{F}$ be a torsion--free coherent sheaf that is $C'^a_{q,k}$ for $(X,H)$ where $q,k$ are integers with $q\geq 0$. Then for any $\zeta\in\textrm{Pic}^0(A)$ and any smooth and irreducible member $Y\in |H+a^*\zeta|$, there exists $\zeta'\in\textrm{Pic}^0(A)$ such that $\mathcal{F}\otimes a^*\zeta'|_Y$ is $C_{q,k}$ for $(Y,H|_Y)$.
\end{lemma}

\noindent\textit{Proof.} It is safe to assume $k=0$. Consider the following part of the long exact sequence
\begin{equation*}
    H^{q+i}(\mathcal{F}(-iH)\otimes a^*\zeta')\to H^{q+i}(\mathcal{F}(-iH)|_Y\otimes a^*\zeta')\to H^{q+i+1}(\mathcal{F}(-(i+1)H- a^*\zeta+ a^*\zeta'))
\end{equation*}
obtained from twisting the restriction sequence by general $\zeta'\in\textrm{Pic}^0(A)$ and passing to cohomology. The statement now follows from semicontinuity.\QEDB

\subsubsection{} We now recall the definition of a {\it strongly generating morphism} as presented in the introduction of \cite{BPS'}.

\begin{definition}
{\bf(Strongly generating morphisms)} Let $a:X\to A$ be a morphism from a smooth projective variety $X$ to an abelian variety $A$. We call the morphism $a$ {\it strongly generating} if the induced map $a^*:\hat{A}=\textrm{Pic}^0(A)\to\textrm{Pic}^0(X)$ is injective.
\end{definition}

Inspired by \cite{Pardini}, given a morphism $a:X\to A$ as above,  
we will work with the {\it covering trick} i.e. we will consider the following base--change diagram 
\begin{equation}\label{trick}
    \begin{tikzcd}
    \Tilde{X}\arrow[r, "\Tilde{\mu}"] \arrow[d, swap, "\Tilde{a}"] & X\arrow[d, "a"]\\
    \Tilde{A}\arrow[r, "\mu"] & A
\end{tikzcd}
\end{equation}
where $\mu:\Tilde{A}\to A$ is an isogeny of abelian varieties. It turns out that if $a$ is strongly generating, then $\Tilde{X}$ is smooth and irreducible (see the proof of \cite[Lemma 2.3]{BPS}). This is precisely the reason we consider strongly generating morphisms. 

Most of the time throughout the article, we will only consider the case when $\Tilde{A}=A$, $\mu=n_A$ is the multiplication by $n$ isogeny for an integer $n\geq 1$, and in this case we denote $\Tilde{X}$ by $X_n$, $\Tilde{a}$ by $a_n$, and $\Tilde{\mu}$ by $\mu_n$. Moreover, in this case, if $a$ is strongly generating then so is $a_n$.  

\smallskip

We show that continuous CM--regularity behaves well with respect to the above covering trick.
\begin{lemma}\label{isogeny}
Let $X$ be a smooth projective variety and let $H$ be a line bundle on $X$. Also, let $a:X\to A$ be a morphism to an abelian variety $A$ that is strongly generating. Let $\mathcal{F}$ be a coherent sheaf on $X$ that is $C'^a_{0,k}$ for $(X,H)$. Let $\mu:\Tilde{A}\to A$ be an isogeny and consider the base--change diagram \eqref{trick}. Then $\Tilde{\mu}^*\mathcal{F}$ is $C'^{\Tilde{a}}_{0,k}$ for $(\Tilde{X},\Tilde{H}:=\Tilde{\mu}^*H)$.
\end{lemma}

\noindent\textit{Proof.} Without loss of generality, we may assume $k=0$. For $\zeta\in\hat{A}$, we have by the projection formula
\begin{equation*}
    H^i(\Tilde{\mu}^*(\mathcal{F}(-iH+a^*\zeta)))=H^i(\mathcal{F}(-iH+a^*\zeta))\oplus\bigoplus\limits_{k=1}^{d-1}H^i(\mathcal{F}(-iH+a^*\zeta+a^*\zeta_k))
\end{equation*}
where $d:=\textrm{deg}(\mu)$, $\mu_*\mathcal{O}_{\Tilde{A}}=\mathcal{O}_{A}\oplus\bigoplus\limits_{k=1}^{d-1}\zeta_k$, and $\zeta_k\in\hat{A}$ for $1\leq k\leq d-1$. The conclusion now follows from semicontinuity and the commutativity of the diagram \eqref{trick}. \QEDB

\section{Continuous CM--regularity and generic vanishing}\label{4.5}

\subsection{The theory of generic vanishing}\label{secgv} Throughout this subsection, for an abelian variety $A$, we denote by $\hat{A}$ the dual abelian variety identified with $\textrm{Pic}^0(A)$.  

\subsubsection{} Let $X$ be a smooth projective variety and let $\textrm{alb}_X:X\to A:=\textrm{Alb}(X)$ be the Albanese map. Consider $\hat{A}=\textrm{Pic}^0(A)\cong \textrm{Pic}^0(X)$ and let $\mathcal{P}$ be a Poincar\'e line bundle on $A\times\hat{A}$. Let $P:=(\textrm{alb}_X\times \textrm{id}_{\hat{A}})^*\mathcal{P}$. Let ${\bf D}(X)$ and ${\bf D}(\hat{A})$ be the bounded derived categories of $\textrm{Coh}(X)$ and $\textrm{Coh}(\hat{A})$. In this situation, we have the following two {\it Fourier--Mukai transform} functors
\begin{align*}
    \begin{split}
        {\bf R}\Phi_{P}:{\bf D}(X)\to{\bf D}(\hat{A}), & \quad {\bf R}\Phi_P(-):={\bf R}{p_{\hat{A}}}_*(p_X^*(-)\otimes P),\\
        {\bf R}\Psi_{P}:{\bf D}(\hat{A})\to{\bf D}(X), & \quad {\bf R}\Psi_P(-):={\bf R}{p_X}_*(p_{\hat{A}}^*(-)\otimes P).
    \end{split}
\end{align*}
A reference of the following definitions can be found for example in \cite[Proposition/Definition 2.1, and Proposition/Definition 2.7]{PaPo3}.

\begin{definition}
{\bf (Generic vanishing and M--regularity)} Let $X$ be a smooth projective variety and let $\mathcal{F}$ be a coherent sheaf on $X$. 
\begin{itemize}
    \item[--] $\mathcal{F}$ is called {\it Generic Vanishing} (abbreviated as {\it GV}) if $\textrm{codim}(V^i(\mathcal{F}))\geq i$ for all integers $i>0$. More generally, for an integer $k\geq 0$, $\mathcal{F}$ is called GV$_{-k}$ if codim$(V^i(\mathcal{F}))\geq i-k$ for all integers $i$.
    \item[--] $\mathcal{F}$ is called {\it Mukai regular} (abbreviated as {\it M--regular}) if $\textrm{codim}(V^i(\mathcal{F}))>i$ for all $i>0$.
\end{itemize} 
\end{definition}

Evidently, we have GV $=$ GV$_0$. The following fundamental theorem is due to Hacon (see \cite{Hac}), Pareschi and Popa (see \cite{PaPo3}).

\begin{theorem}\label{hacon}
(\cite[Theorem 3.7, Corollary 3.11]{PaPo4}) Let $X$ be a smooth projective variety  with $\dim\textrm{Alb(X)}=g$. Let $\mathcal{F}$ be a coherent sheaf on $X$, and let $k\geq 0$ be an integer. Then the following are equivalent:
\begin{itemize}
\item[(1)] $\mathcal{F}$ is GV$_{-k}$,
    \item[(2)] codim Supp$(R^i\Phi_{{P}}\mathcal{F})\geq i-k$ for all integers $i$,
    \item[(3)] for any sufficiently positive ample line bundle $L$ on $\hat{A}$, $H^i(\mathcal{F}\otimes {\bf R}\Psi_{{P}[g]}L^{\otimes -1})=0$ for all integers $i>k$.
\end{itemize}
\end{theorem}

\subsubsection{} Now, let $X$ be an abelian variety of dimension $g$. We recall the notions of M–regularity and Index Theorems with prescribed indices (cf. \cite[the end of Section 1 on p. 5]{PaPo1}).

\begin{definition}\label{defit}
{\bf(IT sheaves)} Let $X$ be an abelian variety and let $\mathcal{F}$ be a coherent sheaf on $X$. The sheaf $\mathcal{F}$ is said to satisfy {\it Index Theorem with index $k$} for some $k\in\mathbb{Z}$ (abbreviated as {\it IT$_{k}$}) if $V^i(\mathcal{F})=\emptyset$ for all $i\neq k$.
\end{definition}

It is clear that on an abelian variety $X$, a coherent sheaf $\mathcal{F}$ satisfies IT$_{0}\implies$ $\mathcal{F}$ is M--regular$\implies$ $\mathcal{F}$ is GV. Also, notice that an ample line bundle on an abelian variety $X$ satisfies IT$_{0}$. In the abelian case, we will denote the Fourier--Mukai transform functors as 
$${\bf R}\hat{\mathcal{S}}:{\bf D}(X)\to {\bf D}(\hat{X}),\quad {\bf R}\mathcal{S}:{\bf D}(\hat{X})\to{\bf D}(X).$$
A fundamental result of Mukai (see \cite[Theorem 2.2]{Muk81}) shows that ${\bf R}\hat{\mathcal{S}}:{\bf D}(X)\to{\bf D}(\hat{X})$ is an equivalence of derived categories, and we have the following inversion formulae
\begin{equation}\label{inversion}
    {\bf R}\hat{\mathcal{S}}\circ{\bf R}\mathcal{S}=(-1)_{\hat{X}}^*[-g],\quad\textrm{ and }\quad {\bf R}\mathcal{S}\circ {\bf R}\hat{\mathcal{S}}=(-1)_{X}^*[-g]
\end{equation}
where $(-1)_{\hat{X}}$ and $(-1)_{X}$ are multiplications by $(-1)$ on $\hat{X}$ and $X$ respectively. 

It turns out that if $\mathcal{F}$ on $X$ is IT$_k$ for some $k\in\mathbb{Z}$, then ${\bf R}\hat{\mathcal{S}}\mathcal{F}=R^k\hat{\mathcal{S}}\mathcal{F}[-k]$ and $R^k\hat{\mathcal{S}}\mathcal{F}$ is a locally free sheaf. In particular, if $L$ is an ample line bundle on $\hat{X}$, then ${\bf R}\mathcal{S}L^{\otimes -1}=R^g\mathcal{S}L^{\otimes -1}[-g]$.

\subsubsection{} A result of Pareschi--Popa on {\it preservation of vanishing} says that on an abelian variety, a tensor product of a GV sheaf $\mathcal{F}$ and an IT$_0$ sheaf $\mathcal{G}$ is IT$_0$ if one of $\mathcal{F}$ and $\mathcal{G}$ is locally free (\cite[Proposition 3.1]{PaPo3}). A variation of their proof yields the following statement, which is more suited for our purposes.

\begin{proposition}\label{preservation}
Let $X$ be an abelian variety of dimension $g$. Let $\mathcal{F}$ and $\mathcal{G}$ be coherent sheaves on $X$ with one of them locally free. Assume $\mathcal{G}$ is IT$_0$ and $m\geq 0$ is an integer. If one of the following holds
\begin{itemize}
    \item[(1)] either $V^i(\mathcal{F})=\emptyset$ for all integers $i\geq m+1$, or
    \item[(2)] $\mathcal{F}$ is GV$_{-m}$,
\end{itemize}
then $V^i(\mathcal{F}\otimes\mathcal{G})=\emptyset$ for all integers $i\geq m+1$.
\end{proposition}

\noindent\textit{Proof.} The proof is identical to that of \cite[Proposition 3.1]{PaPo3}. Let $\zeta\in\textrm{Pic}^0(X)$, and we aim to show that the group $H^i(\mathcal{F}\otimes\mathcal{G}\otimes\zeta)=0$ for all $i\geq m+1$. As $\mathcal{G}\otimes\zeta$ is also IT$_0$, we have ${\bf R}\hat{\mathcal{S}}(\mathcal{G}\otimes\zeta)=R^0\hat{\mathcal{S}}(\mathcal{G}\otimes\zeta)=:\mathcal{N}_{\zeta}$ is locally free. By Mukai's inversion formulae \eqref{inversion}, we have $\mathcal{G}\otimes\zeta={\bf R}\mathcal{S}((-1)^*_{\hat{X}}\mathcal{N}_{\zeta})[g]$. Thus, we deduce that $H^i(X,\mathcal{F}\otimes\mathcal{G}\otimes\zeta)$ is isomorphic to
\begin{equation}\label{exchange}
    H^i(X,\mathcal{F}\otimes{\bf R}\mathcal{S}((-1)^*_{\hat{X}}\mathcal{N}_{\zeta})[g])\cong H^i(\hat{X},{\bf R}\hat{\mathcal{S}}\mathcal{F}\underline{\otimes}(-1)^*_{\hat{X}}\mathcal{N}_{\zeta}[g])\cong H^{g+i}(\hat{X},{\bf R}\hat{\mathcal{S}}\mathcal{F}\underline{\otimes}(-1)^*_{\hat{X}}\mathcal{N}_{\zeta})
\end{equation}
where the first isomorphism in \eqref{exchange} is obtained by an exchange formula of Pareschi--Popa (see \cite[Lemma 2.1]{PaPo4}). To this end, consider the following spectral sequence
\begin{equation*}
    E^{jk}_2:=H^j(\hat{X},R^k\hat{\mathcal{S}}\mathcal{F}\otimes (-1)^*_{\hat{X}}\mathcal{N}_{\zeta})\implies H^{j+k}(\hat{X},{\bf R}\hat{\mathcal{S}}\mathcal{F}\underline{\otimes} (-1)^*_{\hat{X}}\mathcal{N}_{\zeta}).
\end{equation*}
Observe that the $E_2^{jk}$ term above vanishes if $j+k\geq g+m+1$. Indeed, if (1) holds, then $R^k\hat{\mathcal{S}}\mathcal{F}=0$ if $k\geq m+1$, and if (2) holds, then $\dim \textrm{Supp}(R^k\hat{\mathcal{S}}\mathcal{F})\leq g-k+m$. Thus, $E^{jk}_{\infty}=0$ if $j+k\geq g+m+1$, whence $H^{j+k}(\hat{X},{\bf R}\hat{\mathcal{S}}\mathcal{F}\underline{\otimes} (-1)^*_{\hat{X}}\mathcal{N}_{\zeta})=0$ in the same range.\QEDB

\vspace{5pt}

We include a consequence of the above proposition which we will not use anywhere in the sequel.

\begin{corollary}
Let $\mathcal{F}$ and $\mathcal{G}$ be coherent sheaves on an abelian variety $X$ with one of them locally free. If $\mathcal{F}$ is GV$_{-k}$ and $\mathcal{G}$ is GV for some integer $k\geq 0$, then $\mathcal{F}\otimes\mathcal{G}$ is GV$_{-k}$.
\end{corollary}

\noindent\textit{Proof.} This is identical to \cite[Theorem 3.2]{PaPo3}. Let $L$ be a sufficiently positive ample line bundle on $\hat{X}$. Since $\mathcal{G}$ is GV, $\mathcal{G}\otimes R^g\mathcal{S}L^{\otimes -1}$ is IT$_0$ by Theorem \ref{hacon}. Thus, by Proposition \ref{preservation}, $H^i(\mathcal{F}\otimes \mathcal{G}\otimes R^g\mathcal{S}L^{\otimes -1})=0$ for $i\geq k+1$ whence the assertion follows from Theorem \ref{hacon}.\QEDB

\subsubsection{} One of the most essential and useful tools in the study of irregular varieties is the notion of continuous global generation that we define next (\cite[Definition 5.2]{PaPo3}).

\begin{definition}\label{cggprop}
{\bf(Continuous global generation)} Let $X$ be an irregular variety. A sheaf $\mathcal{F}$ on $X$ is {\it continuously globally generated} if for any non--empty open subset $U\subset \textrm{Pic}^0(X)$, the following sum of evaluation maps is surjective:
\begin{equation*}
    \bigoplus\limits_{\zeta\in U}H^0(\mathcal{F}\otimes\zeta)\otimes\zeta^*\to \mathcal{F}.
\end{equation*} 
\end{definition}

It was shown by Pareschi--Popa in \cite[Corollary 5.3]{PaPo3} that an M--regular sheaf is continuously globally generated. Also, another result of them asserts that if $\mathcal{F}$ is a continuously globally generated sheaf and $L$ is a continuously globally generated line bundle on $X$, then $\mathcal{F}\otimes L$ is globally generated (\cite[Proposition 2.12]{PaPo1}). Thus, if $H_1$ and $H_2$ are two ample line bundles on an abelian variety $X$, then $H_1+H_2$ is globally generated. 
It has been pointed out by the referee that alternatively, one could also use \cite[Theorem 1.1]{BS} to get the above statement.

\subsection{Proof of Theorem \texorpdfstring{\hyperref[newgv]{A}}{TEXT}} We now aim to prove Theorem \hyperref[newgv]{A} stated in the introduction. 

\begin{proposition}\label{starit0}
Let $(X,H)$ be a polarized smooth projective variety and let $a:X\to A$ be a morphism to an abelian variety $A$. Assume that $(X,H)$ satisfies $(P_a)$ and let $\mathcal{F}$ be a torsion--free coherent sheaf on $X$ that is $C'^a_{q-1,0}$ for $(X,H)$ for some integer $1\leq q\leq \dim X$. Then $V_a^i(\mathcal{F})=\emptyset$ for all integers $i\geq q$.
\end{proposition}

\noindent\textit{Proof.} The proof is based on induction on $\dim X=:n$ and we divide the proof into the two steps. The following first step verifies the base case which is when $n=1$: 

\smallskip

\noindent\underline{Step 1.} Here we prove the statement for $i=n$. By hypothesis, we know that $H^{n}(\mathcal{F}(-(n-q+1)H)\otimes a^*\zeta)=0$ for some $\zeta\in\textrm{Pic}^0(A)$. Since $(X,H)$ satisfies $(P_a)$, for any given $\zeta'\in\textrm{Pic}^0(A)$, we can choose a smooth and irreducible member $Y\in (n-q+1)H+a^*(\zeta'-\zeta)$. Consider the restriction exact sequence $$0\to \mathcal{F}(-(n-q+1)H)\otimes a^*\zeta)\to \mathcal{F}\otimes a^*\zeta'\to \mathcal{F}\otimes a^*\zeta'|_Y\to 0.$$
Passing to cohomology, we obtain the desired vanishing $H^n(\mathcal{F}\otimes a^*\zeta')=0$.

\smallskip

\noindent\underline{Step 2.} Now we prove the statement by induction and thanks to the previous step, we assume $n\geq 2$. Also, because of Step 1, we assume that $1\leq q\leq i\leq n-1$. Let $\zeta\in\textrm{Pic}^0(A)$ and we want to show that $H^i(\mathcal{F}\otimes a^*\zeta)=0$. Our proof is inspired by the proof of \cite[Lemma 3.3]{LS}. We know by Observation \ref{obs} that there exists $\zeta'\in\textrm{Pic}^0(A)$ such that $\mathcal{F}\otimes a^*\zeta'$ is $C_{q-1,0}$ for $(X,H)$. Observe that by Lemma \ref{partial}, $H^i(\mathcal{F}\otimes a^*\zeta'\otimes (-jH))=0$ for all integers $j\leq i-q+1$. To this end, consider the restriction exact sequence
\begin{equation}\label{restrseq}
    0\to\mathcal{F}\otimes a^*\zeta\to\mathcal{F}\otimes (H+a^*\zeta')\to\mathcal{F}\otimes (H+ a^*\zeta')|_{Y}\to 0
\end{equation} where $Y\in|H+a^*\zeta'-a^*\zeta|$ is a smooth and irreducible member (exists by Bertini thanks to $(P_a)$). Passing to the cohomology of \eqref{restrseq}, we deduce that it is enough to prove that the restriction map $$H^{i-1}(\mathcal{F}\otimes(H+ a^*\zeta'))\to H^{i-1}(\mathcal{F}\otimes(H+a^*\zeta')|_Y)$$ surjects since $H^i(\mathcal{F}\otimes (H+a^*\zeta'))=0$. On the other hand, choose a general smooth and irreducible member $Z\in |H+a^*\zeta-a^*\zeta'|$ such that $\mathcal{F}|_Z$ is torsion--free (such a section exists thanks again to the fact that $(X,H)$ satisfies property $(P_a)$). Then we know that $H^{i-1}(\mathcal{F}\otimes (H+a^*\zeta'))\to H^{i-1}(\mathcal{F}\otimes (H+a^*\zeta')|_{Y+Z})$ surjects by Lemma \ref{partial} as $H^i(\mathcal{F}\otimes a^*\zeta'(-H))=0$. Consequently, it is enough to show that the map $$H^{i-1}(\mathcal{F}\otimes(H+ a^*\zeta')|_{Y+Z})\to H^{i-1}(\mathcal{F}\otimes(H+ a^*\zeta')|_Y)$$ surjects. Consider the following commutative diagram with exact rows and exact left column
\small 
\[
\begin{tikzcd}
    & 0\arrow[d]&&&\\
    0 \arrow{r} & \mathcal{O}_X(-Y-Z) \arrow[r]\arrow[d] & \mathcal{O}_X \arrow{r} \arrow[d, equal] & \mathcal{O}_{Y+Z}\arrow[r]\arrow[d] & 0\\
   0 \arrow{r} & \mathcal{O}_X(-Y) \arrow{r}\arrow[d] & \mathcal{O}_X \arrow{r} & \mathcal{O}_Y\arrow[r] & 0\\
   & \mathcal{O}_Z(-Y)\arrow[d] &&&\\
   & 0 &&&
\end{tikzcd} \]
\normalsize
which by snake lemma yields the following short exact sequence
\begin{equation*}
    0\to\mathcal{O}_Z(-Y)\to\mathcal{O}_{Y+Z}\to\mathcal{O}_Y\to 0.
\end{equation*}
Twisting the above by $\mathcal{F}\otimes (H+a^*\zeta')$ and passing to cohomology, we deduce that it is enough to prove the following vanishing $H^i(\mathcal{F}\otimes (H+ a^*\zeta')\otimes (-(H+a^*\zeta'-a^*\zeta))|_Z)=H^i(\mathcal{F}\otimes a^*\zeta|_Z)=0.$
But $\mathcal{F}|_Z$ is torsion--free and $C'^{a|_Z}_{q-1,0}$ for $(Z,H|_Z)$ by Lemma \ref{contres}, and $(Z,H|_Z)$ satisfies $(P_{a|_Z})$ by Observation \ref{obs2}. Thus we are done by the induction hypothesis.\QEDB

\begin{theorem}\label{unk}
Let $(X,H)$ be a polarized smooth projective variety and let $a:X\to A$ be a morphism to an abelian variety $A$. Assume that $(X,H)$ satisfies $(P_a)$. Let $\mathcal{F}$ be a torsion--free coherent sheaf on $X$ that is $C'^a_{q-1,0}$ for $(X,H)$ for some integer $1\leq q\leq \dim X$. Then for any integer $0\leq i\leq \dim X-q$, $V_a^{q+i}(\mathcal{F}(-t H))=\emptyset$ for all integers $0\leq t\leq i$.
\end{theorem}

\noindent\textit{Proof.} Clearly the statement holds by Proposition \ref{starit0} if $t=0$, in particular it holds if $\dim X=:n=1$. Consider the restriction sequence $$0\to\mathcal{F}(-t H+a^*\zeta)\to\mathcal{F}(-(t-1)H+a^*\zeta)\to\mathcal{F}(-(t-1)H+a^*\zeta)|_Y)\to 0$$ where $Y\in|H|$ is a general smooth and irreducible member such that $\mathcal{F}|_Y$ is torsion--free. The cohomology sequence of the above, and an easy induction on $n$ and $t$ finishes the proof (thanks to Lemma \ref{contres} and Observation \ref{obs2}).\QEDB

\begin{corollary}\label{newgvcor}
Let $(X,H)$ be a polarized smooth projective variety and let $a:X\to A$ be a morphism to an abelian variety $A$. Assume one of the following two conditions holds:
\begin{itemize}
    \item[(1)] $(X,H)$ satisfies property $(P_a)$, or
    \item[(2)] the morphism $a$ is strongly generating, and there exists an isogeny $\mu:\Tilde{A}\to A$ such that $(\Tilde{X},\Tilde{\mu}^*H)$ satisfies $(P_{\Tilde{a}})$, where $\Tilde{X}$, $\Tilde{a}$, and $\Tilde{\mu}$ are as in \eqref{trick}. 
\end{itemize}
Let $\mathcal{F}$ be a torsion--free coherent sheaf on $X$ that is $C'^a_{0,k}$ for some integer $0\leq k\leq \dim X$. Then 
\begin{itemize}
    \item[(A)] $V^i_a(\mathcal{F})=\emptyset$ for all integers $i\geq k+1$,
    \item[(B)] if $k\neq 0$, then $\textrm{codim}(V^k_a(\mathcal{F}))\geq 1$.
\end{itemize}
\end{corollary}

\noindent\textit{Proof.} The assertion follows immediately from Theorem \ref{unk} if (1) holds. Now assume (2) holds. The assertion (B) is obvious, so it is enough to show (A), that is $V^i(\mathcal{F})=\emptyset$ for $i\geq k+1$. By Lemma \ref{isogeny}, we know that $\Tilde{\mu}^*\mathcal{F}$ is $C'^{\Tilde{a}}_{0,k}$ for $(\Tilde{X},\Tilde{\mu}^*H)$. Also, $\Tilde{\mu}^*\mathcal{F}$ is torsion--free and coherent. Since $(\Tilde{X},\Tilde{\mu}^*H)$ satisfies $(P_{\Tilde{a}})$, using Theorem \ref{unk} we conclude that $V^i_{\Tilde{a}}(\Tilde{\mu}^*\mathcal{F})=\emptyset$ for $i\geq k+1$. It follows from the projection formula that $\mu^*V^i_a(\mathcal{F})\subseteq V^i_{\Tilde{a}}(\Tilde{\mu}^*\mathcal{F})$ for all $i$ where $\mu^*:\textrm{Pic}^0(A)\to\textrm{Pic}^0(A)$ is the induced map, whence the assertion follows.\QEDB

\vspace{5pt}

We are now ready to provide the

\vspace{5pt}

\noindent\textit{Proof of Theorem \hyperref[newgv]{A}.} This is an immediate consequence of Corollary \ref{newgvcor} (2). Indeed, set $a:=\textrm{alb}_X$, $A:=\textrm{Alb}(X)$, and let $n_A:A\to A$ be the multiplication by $n$ isogeny. Then $a$ is strongly generating. Observe that $\mu_n^*H=\mu_n^*H_1+a_n^*(n_A^*H_2)$, and consequently $(X_n,\mu_n^*H)$ satisfies $(P_{a_n})$ for $n\geq 2$. \QEDB

\subsection{A few variants}\label{secvar} In this subsection, we prove several other variants of Theorem \hyperref[newgv]{A}.

\subsubsection{} We use the following corollary to show that the moduli space of Gieseker--stable sheaves on abelian surfaces answers Question \ref{newq} in the affirmative.

\begin{corollary}\label{var1}
Let $X$ be a smooth projective variety, and let $\textrm{alb}_X:X\to A:=\textrm{Alb}(X)$ be the Albanese map of $X$. Assume there is an isogeny $\mu:\Tilde{A}\to A$ such that the following two conditions hold: 
\begin{itemize}
    \item[(1)] $X\times_{\textrm{Alb}(X)}\Tilde{A}\cong Y\times \Tilde{A}$ for a regular smooth projective variety $Y$, and
    \item[(2)] the induced map $X\times_{\textrm{Alb}(X)}\Tilde{A}\to \Tilde{A}$ is the projection $\textrm{pr}_{\Tilde{A}}$ under the identification in (1).
\end{itemize} 
Let $\mathcal{O}_X(1)$ is an ample and globally generated line bundle on $X$. Let $\mathcal{F}$ be a torsion--free coherent sheaf on $X$ that is continuously $k$--regular for $(X,\mathcal{O}_X(1))$ where $0\leq k\leq \dim X$. Then 
\begin{itemize}
    \item[(A)] $V^i(\mathcal{F})=\emptyset$ for $i\geq k+1$.
    \item[(B)] If $k\neq 0$, then $\textrm{codim}(V^k(\mathcal{F}))\geq 1$.
\end{itemize}
In particular, Question \ref{newq} has an affirmative answer for $X$.
\end{corollary}

\noindent\textit{Proof.} Since (B) is obvious, we only need to show (A). We have the following base--change diagram 
\[
\begin{tikzcd}
    \Tilde{X}:=Y\times \Tilde{A}\arrow[r, "\Tilde{\mu}"]\arrow[d, swap, "\Tilde{a}=\textrm{pr}_{\Tilde{A}}"] & X\arrow[d,"a:=\textrm{alb}_X"]\\
    \Tilde{A}\arrow[r, "\mu"] & A
\end{tikzcd}
\]
Set $H=\mathcal{O}_X(1)$. Notice that since $Y$ is regular, $\Tilde{a}$ is also the Albanese map of $\Tilde{X}$. It is easy to verify that it is enough to show that $V^i(\Tilde{\mu}^*(\mathcal{F}))=\emptyset$ for $i\geq k+1$. To this end, we apply Lemma \ref{isogeny} to deduce that $\Tilde{\mu}^*(\mathcal{F})$ is a torsion--free coherent sheaf that is continuously $k$--regular. Now, $\Tilde{\mu}^*(H)=H_1\boxtimes H_2$ where $H_1$ and $H_2$ are ample and globally generated line bundles on $Y$ and $\Tilde{A}$ respectively. Consequently the assertion follows from Theorem \hyperref[newgv]{A}.
\QEDB

\begin{example}
Let $A$ be an abelian surface and let $v\in H^{\textrm{ev}}(A,\mathbb{Z})$ be a primitive Mukai vector satisfying $v>0$ and $\langle v,v\rangle\geq 6$ (see \cite{Yos} for details). Let $L$ be a very general ample divisor on $A$, and let $M_L(v)$ be the moduli space of Gieseker--stable sheaves on $A$ with respect to $L$ with Mukai vector $v$. By the results of {\it loc. cit.}, we know that $M_L(v)$ is a smooth projective variety, and $\textrm{Alb}(M_L(v))=A\times\hat{A}$. Moreover, if we set $n(v)=\frac{1}{2}\langle v,v\rangle$, and $a:=\textrm{alb}_{M_L(v)}$, then we have the base--change diagram (see {\it loc. cit.} (4.10), (4.11))
\[
\begin{tikzcd}
    K_L(v)\times A\times \hat{A}\arrow[r, "\mu_{n(v)}"]\arrow[d, swap, "a_{n(v)}"] & M_L(v)\arrow[d,"a"]\\
    A\times\hat{A}\arrow[r, "n(v)_{A\times\hat{A}}"] & A\times\hat{A}
\end{tikzcd}
\]
where $K_L(v)$ is a regular smooth projective variety (it is a hyperk\"ahler manifold deformation equivalent to generalized Kummer variety), and $a_{n(v)}=\textrm{pr}_{A\times\hat{A}}$ is the projection. Thus Question \ref{newq} has an affirmative answer for $X=M_L(v)$ by the previous corollary. 
\end{example}

We have the following example when the polarization is a sufficiently positive adjoint linear series:

\begin{example}
Let $X$ be a smooth projective variety and let $H\equiv K_X+sL$ where $s\geq \dim X+1$ and $L$ is an ample and globally generated line bundle on $X$. By Corollary \ref{newgvcor} (1), the answer of Question \ref{newq} is affirmative for the pair $(X,H)$ whenever $H$ is ample.
\end{example}

\subsubsection{} The following result is a variant of Theorem \hyperref[newgv]{A} where we assume $H_1=K_X+Q$ for a nef line bundle $Q$, but require that $\textrm{alb}_X:X\to\textrm{Alb}(X)$ is finite onto its image.

\begin{theorem}\label{var2}
Let $(X,H)$ be a polarized smooth projective variety. Assume that the Albanese map $\textrm{alb}_X:X\to\textrm{Alb}(X)$ is finite onto its image. Further assume that there exists a nef line bundle $Q$ on $X$ and an ample line bundle $H_2$ on $\textrm{Alb}(X)$ such that $H=K_X+Q+\textrm{alb}_X^*H_2$. Let $\mathcal{F}$ be a torsion--free coherent sheaf on $X$ that is continuously $k$--regular for $(X,H)$ for some integer $0\leq k\leq \dim X$. Then the following statements hold.
\begin{itemize}
    \item[(1)] $V^i(\mathcal{F})=\emptyset$ for $i\geq k+1$.
    \item[(2)] If $k\neq 0$, then $\textrm{codim}(V^k(\mathcal{F}))\geq 1$.
\end{itemize} 
In particular, the answer of Question \ref{newq} is affirmative for the pair $(X,H)$.
\end{theorem}

\noindent\textit{Proof.} This is also an immediate consequence of Corollary \ref{newgvcor} (2). Indeed, set $a:=\textrm{alb}_X$ and $A:=\textrm{Alb}(X)$, and notice that $\mu_n^*H\equiv K_{X_n}+\mu_n^*(Q)+a_n^*(n^2H_2)$ satisfies $(P_{a_n})$ for $n\gg 0$.\QEDB

\subsubsection{} Let $(X,H)$ be a polarized smooth projective variety, and let $\mathcal{F}$ be a torsion--free coherent sheaf on $X$ that is continuously $k$--regular as in the statements of Theorems \hyperref[newgv]{A}, \ref{var2} and Corollary \ref{newgvcor}. Then in particular we have showed that $\mathcal{F}$ is GV$_{-(k-1)}$ whence by Theorem \ref{hacon} $H^i(\mathcal{F}\otimes {\bf R}\Psi_{P[g]}L^{\otimes -1})=0$ for all $i\geq k$ and for any sufficiently positive ample line bundle $L$ on $\hat{A}$. We now show that in some cases the vanishing in fact holds for {\it any} ample line bundle $L$ on $\hat{A}$. First we need the following result.

\begin{corollary}\label{earlier}
Let $a:X\to A$ be a morphism from a smooth projective variety $X$ to an abelian variety $A$ that is finite onto its image. Let $H$ be an ample line bundle on $X$ that satisfies one of the following:
\begin{itemize}
    \item[(1)] $(X,H)$ satisfies property $(P_a)$; or 
    \item[(2)] $a$ is strongly generating, and moreover $H:=H_1+a^*H_2$ where $H_1$ is a line bundle on $X$ and $H_2$ is an ample line bundle on $A$ satisfying one of the following:
    \begin{itemize}
        \item[(i)] $H_1$ is globally generated; or
        \item[(ii)] $H_1=K_X+Q$ for a nef line bundle $Q$.
    \end{itemize}
\end{itemize}
Let $\mathcal{F}$ be a torsion--free coherent sheaf on $X$ that is $C'^a_{0,k}$ for $(X,H)$ for some integer $1\leq k\leq \dim X$.
Then $H^i(\mathcal{F}\otimes a^*L)=0$ for any ample line bundle $L$ on $A$ and for any integer $i\geq k$. 
\end{corollary}

\noindent\textit{Proof.} It follows from the proofs of Theorems \hyperref[newgv]{A}, \ref{var2} and Corollary \ref{newgvcor} that $\mathcal{F}$ is GV$_{-(k-1)}$ with respect to $a$ (this means that codim$V^i_a(\mathcal{F})\geq i-k+1$ for all $i$). Consequently $H^i(\mathcal{F}\otimes a^*L)=0$ for $i\geq k$ and any ample line bundle $L$ on $A$ follows from the projection formula and Proposition \ref{preservation}.\QEDB

\vspace{5pt}

The above vanishing gives the following

\begin{corollary}
    Let $X$ be a smooth projective variety and let $a:=\textrm{alb}_X:X\to A:=\textrm{Alb}(X)$ be its Albanese morphism. Assume $a$ is finite onto its image and let $H$ be an ample line bundle on $X$ such that $(X,H)$ satisfies $(P_a)$. Let $\mathcal{F}$ be a torsion--free coherent sheaf on $X$ that is $C'^a_{0,k}$ for $(X,H)$ for some integer $1\leq k\leq \dim X$. Then $H^i(\mathcal{F}\otimes {\bf R}\Psi_{P[g]}L^{\otimes -1})=0$ for all $i\geq k$ and for any ample line bundle $L$ on $\hat{A}$ where $g=q(X)$.
\end{corollary}

\noindent\textit{Proof.} Fix an ample line bundle $L$ on $\hat{A}$. We recall from the proof of \cite[Theorem B]{PaPo4} that ${\bf R}\Psi_{P[g]}L^{\otimes -1}= \textrm{alb}_X^*(R^0\mathcal{S}L)^*$. Denoting translations by an element $y\in\hat{A}$ by $t_y:\hat{A}\to\hat{A}$, we obtain an isogeny $\varphi_L:\hat{A}\to A$ by sending $y\in \hat{A}$ to $t_y^*L\otimes L^{\otimes -1}$. To this end, consider the following base--change diagram.
\begin{equation}\label{hacondiagram}
\begin{tikzcd}
    \hat{X}\arrow[r, "\hat{\varphi}"]\arrow[d, swap, "\hat{a}"] & X\arrow[d, "\textrm{alb}_X"] \\
    \hat{A}\arrow[r, "\varphi_L"] & A
\end{tikzcd}
\end{equation}
On the other hand, \cite[Proposition 3.11]{Muk81} shows that $\varphi_L^*(R^0\mathcal{S}L)^*=H^0(L)\otimes L$. Consequently, we deduce using the projection formula that there is an injection
\begin{equation*}
    H^i(\mathcal{F}\otimes {\bf R}\Psi_{P[g]}L^{\otimes -1})=H^i(\mathcal{F}\otimes \textrm{alb}_X^*(R^0\mathcal{S}L)^*)\hookrightarrow H^i(\hat{\varphi}^*(\mathcal{F}\otimes \textrm{alb}_X^*(R^0\mathcal{S}L)^*)=H^0(L)\otimes H^i((\hat{\varphi}^*\mathcal{F})\otimes \hat{a}^*L).
\end{equation*}
But $\hat{\varphi}^*\mathcal{F}$ is a torsion--free coherent sheaf that is $C'^{\hat{a}}_{0,k}$ for $(\hat{X},\hat{\varphi}^*H)$ by Lemma \ref{isogeny}. Observe that 
$\hat{\varphi}^*H$ satisfies $(P_{\hat{a}})$, as the induced map $\textrm{Pic}^0(A)\to \textrm{Pic}^0(\hat{A})$ is also an isogeny.
Consequently, $H^i((\hat{\varphi}^*\mathcal{F})\otimes \hat{a}^*L)=0$ for all integers $i\geq k$ by Corollary \ref{earlier} and the assertion follows.\QEDB

\section{Continuous \texorpdfstring{$\mathbb{Q}$}{TEXT} CM--regularity on abelian varieties}\label{4.7}
Throughout this section, $X$ is an abelian variety of dimension $g$, $\mathcal{E}$ is a vector bundle on $X$.
\subsection{Continuous \texorpdfstring{$\mathbb{Q}$}{TEXT} CM--regularity for vector bundles} In view of the recent development of the {\it cohomological rank function} by Jiang--Pareschi (see \cite{JP}) that was motivated by the {\it continuous rank function} introduced and studied in \cite{Bar}, \cite{BPS'}, it is natural to extend the notion of continuous CM--regularity to an $\mathbb{R}$--valued regularity function on $\mathbb{Q}$--twisted vector bundles.  

\begin{definition}
{\bf (Cohomological rank function)} Define $h^i_{\textrm{gen}}(\mathcal{E})$ for $i\in\mathbb{N}$ as the dimension of $H^i(\mathcal{E}\otimes\zeta)$ for general $\zeta\in\textrm{Pic}^0(X)$. Given a polarization $\underline{l}\in N^1(X)$, and $x=a/b\in\mathbb{Q}$ with $a,b\in\mathbb{Z}$ and $b>0$, following \cite[Definition 2.1]{JP}, define the {\it cohomological rank function} $h^i_{\mathcal{E}}(x\underline{l})=b^{-2g}h^i_{\textrm{gen}}(b_X^*\mathcal{E}\otimes L^{\otimes ab})$ where $L$ is a line bundle representing $\underline{l}$.
\end{definition}

We note that if $L$ is an ample line bundle, then $\mathcal{E}$ being continuously $k$--regular for $(X,L)$ for $k\in\mathbb{Z}$ is equivalent to the condition $h^i_{\mathcal{E}}((k-i)\underline{l})=0$ for all integers $i\geq 1$.

\smallskip

Let $\underline{l}\in N^1(X)$ be a polarization. Assume for some $y=a/b\in\mathbb{Q}$ with $a,b\in \mathbb{Z}$ and $b>0$, $h^i_{\mathcal{E}}((y-i)\underline{l})=0$ for all integers $i\geq 1$. This means that $h^i_{\textrm{gen}}(b_X^*\mathcal{E}\otimes L^{\otimes (a-ib)b})=0$ for all integers $i\geq 1$. This is equivalent to the condition that
$b_X^*\mathcal{E}\otimes L^{\otimes ab}$ is continuously $0$--regular for $(X,L^{\otimes b^2})$.
We claim that for any integers $c,d>0$, $h^i_{\mathcal{E}}((\frac{a}{b}+\frac{c}{d}-i)\underline{l})=0$ for all integers $i\geq 1$. To see this, we need to show that $h^i_{\textrm{gen}}((bd)_X^*\mathcal{E}\otimes L^{\otimes (ad+bc-ibd)bd})=0$ for $i\geq 1$, which is equivalent to 
$(bd)_X^*\mathcal{E}\otimes L^{\otimes (ad+bc)bd}$ being continuously $0$--regular for $(X,L^{\otimes (bd)^2})$.
But the given condition implies that $(bd)_X^*\mathcal{E}\otimes L^{\otimes abd^2}$ is continuously $0$--regular for $(X,L^{\otimes (bd)^2})$, and consequently the required continuous regularity follows by Corollary \ref{earlier}. Thus we define the following $$\mathbb{Q}\textrm{-reg}_{\underline{l}}(\mathcal{E}):=\textrm{inf}\left\{y\in\mathbb{Q}\mid h^i_{\mathcal{E}}((y-i)\underline{l})=0\textrm{ for all integers $i\geq 1$}\right\}.$$
\begin{example}
{\it (Continuous $\mathbb{Q}$ CM--regularity for Verlinde bundles)} We compute an example of $\mathbb{Q}$-$\textrm{reg}_{\underline{l}}(\mathcal{E})$ when $X:=J(C)$ is the Jacobian of a smooth projective curve $C$ of genus $g\geq 1$, $\underline{l}=s\underline{\theta}$ is the class of $s\Theta$ where $\Theta$ is a symmetric theta--divisor on $X$, $s\ge 2$ is an integer, and $\mathcal{E}:=\mathbb{E}_{r,k}$ is the Verlinde bundle. We recall the definition of these bundles. For a pair of positive integers $(r,k)$, let $U_C(r,0)$ be the moduli space of semistable bundles of rank $r$ and degree $0$ on $C$, and let $\textrm{det}:U_r:=U_C(r,0)\to\hat{X}=X$ be the determinant map. The {\it Verlinde bundle} associated to $(r,k)$ is by definition $\mathbb{E}_{r,k}:=\textrm{det}_*\mathcal{O}_{U_r}(k\Tilde{\Theta})$ where $\Tilde{\Theta}$ is the generalized theta--characteristic of $C$. For details about these bundles, we refer to \cite{Pop} and \cite{Opr}. 

K\"uronya and Mustopa showed in \cite[Proposition 3.2]{KM} that $\textrm{reg}_{\mathcal{O}_{X}(s\Theta)}^{\textrm{cont}}(\mathbb{E}_{r,k})=\lceil g-\frac{k}{rs}\rceil$ if $2\nmid \textrm{gcd}(r,k)$. We claim that $\mathbb{Q}\textrm{-reg}_{s\underline{\theta}}(\mathbb{E}_{r,k})=g-\frac{k}{rs}$ if $2\nmid \textrm{gcd}(r,k)$. 

We show this by following their proof. First of all, by the proof of \cite[Proposition 3.2]{KM}, one can write $\mathbb{E}_{r,k}=\bigoplus\mathbb{W}_{a,b}\otimes\zeta_i$ for $\zeta_i\in\hat{X}$ where $a=r/\textrm{gcd}(r,k)$, $b=k/\textrm{gcd}(r,k)$, and $\mathbb{W}_{a,b}$ is semihomogeneious i.e., for any $x\in X$, there exists $\xi\in\hat{X}$ such that $t_x^*\mathbb{W}_{a,b}=\mathbb{W}_{a,b}\otimes\xi$. We also know that $\textrm{rank}(\mathbb{W}_{a,b})=a^g$, $\textrm{det}(\mathbb{W}_{a,b})=\mathcal{O}_X(a^{g-1}b\Theta)$. Observe that we have the following equality: 
$$\mathbb{Q}\textrm{-reg}_{s\underline{\theta}}(\mathbb{E}_{r,k})=\textrm{inf}\left\{y\in\mathbb{Q}\mid h^i_{\mathbb{E}_{r,k}}\left(\left(y-i\right)s\underline{\theta}\right)=0\,\forall\, i\geq 1\right\}=\textrm{inf}\left\{y\in\mathbb{Q}\mid h^i_{\mathbb{W}_{a,b}}\left(\left(y-i\right)s\underline{\theta}\right)=0\,\forall\, i\geq 1\right\}$$
$$=\textrm{inf}\left\{\frac{r'}{s'}\in\mathbb{Q}\bigg\vert \textrm{ $s'>0$ and } (s')_X^*\mathbb{W}_{a,b}\otimes (s\Theta)^{\otimes r's'}\textrm{ is continuously $0$--regular for $(X,s'^2s\Theta)$}\right\}.$$
Now $(s')_X^*\mathbb{W}_{a,b}\otimes (s\Theta)^{\otimes r's'}$ is semihomogeneous and $c_1((s')_X^*\mathbb{W}_{a,b}\otimes (s\Theta)^{\otimes r's'})\in N^1(X)$ is $\mathbb{Q}$--proportional to $s'^2s\underline{\theta}$. Thus, by \cite[Proposition 2.8]{KM}, $\mathbb{Q}\textrm{-reg}_{s\underline{\theta}}(\mathbb{E}_{r,k})\leq \frac{r'}{s'}$ if and only if $(s')_X^*\mathbb{W}_{a,b}\otimes \Theta^{\otimes r's's}(-gs'^2s\Theta)$ is nef, which in turn holds if and only if $a^{g-1}(s'^2b+a(r's's-gs'^2s))\geq 0$ by \cite[Proposition 2.7]{KM}. A simple computation shows $a^{g-1}(s'^2b+a(r's's-gs'^2s))\geq 0\iff \frac{r'}{s'}\geq g-\frac{k}{rs}$. 
\end{example}

\subsection{Continuous \texorpdfstring{$\mathbb{Q}$}{TEXT} CM--regularity for \texorpdfstring{$\mathbb{Q}$}{TEXT}--twisted bundles} Given a class $\underline{n}\in\textrm{Pic}(X)/\textrm{Pic}^0(X)=:N^1(X)$, following \cite[Chapter 6]{Laz'}, define the $\mathbb{Q}$--twisted bundle $\mathcal{E}\langle x\underline{n}\rangle$ for $x\in\mathbb{Q}$ as the equivalence class of pairs $(\mathcal{E},x\underline{n})$ with respect to the equivalence relation generated by declaring $(\mathcal{E}\otimes M^{\otimes e},y\underline{n})\sim (\mathcal{E}, e\underline{m}+x\underline{n})$ where $M\in\textrm{Pic}(X)$, $\underline{m}$ is its class in $N^1(X)$, $e\in\mathbb{Z}$ and $y\in\mathbb{Q}$. Now, given a polarization $\underline{l}\in N^1(X)$, we define the {\it continuous $\mathbb{Q}$ CM--regularity} of  $\mathcal{E}\langle x\underline{n}\rangle$ as follows
$$\mathbb{Q}\textrm{-reg}_{\underline{l}}(\mathcal{E}\langle x\underline{n}\rangle):=\mathbb{Q}\textrm{-reg}_{b^2\underline{l}}(b_X^*\mathcal{E}\otimes N^{\otimes ab})$$ where $x=a/b$, $a,b\in\mathbb{Z}$ and $b>0$. It is a formal verification that this quantity is well--defined. We only show that it does not depend on the representation of the $\mathbb{Q}$--twist. The fact that it does not depend on its representation as a $\mathbb{Q}$--twisted bundle under the equivalence described above can also be easily checked. 

Let $\frac{a}{b}\underline{n_1}=\frac{c}{d}\underline{n_2}\in N^1(X)_{\mathbb{Q}}$ with $a,b,c,d\in\mathbb{Z}$, $b,d>0$. Set $\underline{n'}:=ad \underline{n_1}=bc\underline{n_2}\in N^1(X)$. Then we have  
\begin{equation*}
h^i_{b_X^*\mathcal{E}\otimes N_1^{\otimes ab}}((r/s-i)b^2\underline{l})=0 \,\forall\, i\geq 1\iff h^i_{\textrm{gen}}((bs)_X^*\mathcal{E}\otimes N_1^{\otimes s^2ab}\otimes L^{\otimes (r-is)b^2s})=0\,\forall\,i\geq 1\iff
\end{equation*}
\begin{equation*}
    h^i_{\textrm{gen}}((bsd)_X^*\mathcal{E}\otimes N_1^{\otimes d^2s^2ab}\otimes L^{\otimes (r-is)b^2d^2s})=0\,\forall\,i\geq 1\iff h^i_{\textrm{gen}}((bsd)_X^*\mathcal{E}\otimes N'^{\otimes bds^2}\otimes L^{\otimes (r-is)b^2d^2s})=0\,\forall\,i\geq 1.
\end{equation*}
Similar computation as above shows that $h^i_{d_X^*\mathcal{E}\otimes N_2^{\otimes cd}}((r/s-i)d^2\underline{l})=0$ for all $i\geq 1$ is equivalent to
\begin{equation*}
    h^i_{\textrm{gen}}((sd)_X^*\mathcal{E}\otimes N_2^{\otimes s^2cd}\otimes L^{\otimes (r-is)d^2s})=0\,\forall\,i\geq 1\iff h^i_{\textrm{gen}}((bsd)_X^*\mathcal{E}\otimes N_2^{\otimes b^2s^2cd}\otimes L^{\otimes (r-is)b^2d^2s})=0\,\forall\,i\geq 1
\end{equation*}
which is equivalent to the condition $h^i_{\textrm{gen}}((bsd)_X^*\mathcal{E}\otimes N'^{\otimes bds^2}\otimes L^{\otimes (r-is)b^2d^2s})=0\,\forall\,i\geq 1$.

\subsection{Generic vanishing for \texorpdfstring{$\mathbb{Q}$}{TEXT}--twisted bundles} We mention an immediate corollary of Theorem \hyperref[newgv]{A}. In \cite{JP}, Jiang and Pareschi extended the definitions of GV, M--regular and IT$_0$ sheaves to the $\mathbb{Q}$--twisted cases. In particular, according to their definitions, for an ample class $\underline{l}\in N^1(X)$, a $\mathbb{Q}$--twisted vector bundle $\mathcal{E}\langle x\underline{l}\rangle$ for $x=\frac{a}{b}$ with $b>0$ is GV, M--regular, or IT$_0$ if so is $b_X^*\mathcal{E}\otimes L^{\otimes ab}$. 

\begin{corollary}\label{qtwist}
Let $X$ be an abelian variety. Let $\underline{l},\underline{l'}\in N^1(X)$ be ample classes and let $\mathcal{E}$ be a vector bundle on $X$. If $\mathbb{Q}\textrm{-reg}_{\underline{l}}(\mathcal{E}\langle x\underline{l'}\rangle)<1$, then $\mathcal{E}\langle x\underline{l'}\rangle$ is IT$_0$, and if $\mathbb{Q}\textrm{-reg}_{\underline{l}}(\mathcal{E}\langle x\underline{l'}\rangle)=1$, then $\mathcal{E}\langle x\underline{l'}\rangle$ is GV.
\end{corollary}

\noindent\textit{Proof.} Let $x=\frac{r}{s}$ with $s>0$. First, if $\mathbb{Q}\textrm{-reg}_{\underline{l}}(\mathcal{E}\langle x\underline{l'}\rangle)<1$, then we can find $a,b\in\mathbb{Z}$ with $b>0$ such that we have $\mathbb{Q}\textrm{-reg}_{\underline{l}}(\mathcal{E}\langle x\underline{l'}\rangle)<\frac{a}{b}<1$. This means that $h^i_{s_X^*\mathcal{E}\otimes L'^{\otimes ab}}((\frac{a}{b}-i)s^2\underline{l})=0$ for all $i\geq 1$. But this means that the bundle $(bs)_X^*\mathcal{E}\otimes L'^{\otimes b^2rs}\otimes L^{\otimes abs^2}$ is continuously $0$--regular for $(X,L^{\otimes b^2s^2})$. Consequently, by Theorem \hyperref[newgv]{A}, we deduce that $(bs)_X^*\mathcal{E}\otimes L'^{\otimes b^2rs}\otimes ((abs^2-b^2s^2)L)$ is GV, whence by preservation of vanishing, we conclude that $(bs)_X^*\mathcal{E}\otimes L'^{\otimes b^2rs}$ is IT$_0$. Thus, $s_X^*\mathcal{E}\otimes L'^{\otimes rs}$ is IT$_0$, and the conclusion follows. Finally, if $\mathbb{Q}\textrm{-reg}_{\underline{l}}(\mathcal{E}\langle x\underline{l'}\rangle)=1$, then similarly we see that $s_X^*\mathcal{E}\otimes L'^{\otimes rs}$ is GV, and the conclusion follows.\QEDB

\subsection{Proof of Theorem \texorpdfstring{\hyperref[qregthm]{B}}{TEXT}} The proof of Theorem \hyperref[qregthm]{B} is based on the following two results.

\begin{lemma}\label{cont1}
Let $X$ be an abelian variety of dimension $g$ and let $\mathcal{E}$ be a vector bundle on $X$. Let $a,b\in \mathbb{Z}$ with $b>0$, and set $x=a/b\in\mathbb{Q}$. If $\underline{n},\underline{l}\in N^1(X)$ with $\underline{l}$ ample, and $\delta=c/d\in\mathbb{Q}$ with $c,d\in\mathbb{Z}$, $d>0$, then $$\mathbb{Q}\textrm{-reg}_{\underline{l}}(\mathcal{E}\langle x\underline{n}+\delta\underline{l}\rangle)=\mathbb{Q}\textrm{-reg}_{\underline{l}}(\mathcal{E}\langle x\underline{n}\rangle)-\delta.$$
\end{lemma}

\noindent\textit{Proof.} Let $r,s\in\mathbb{Z}$ with $s>0$. It is enough to show the following equivalence 
\begin{equation}\label{q0}
    h^i_{(bd)_X^*\mathcal{E}\otimes N^{\otimes abd^2}\otimes L^{\otimes b^2cd}}((r/s-i)b^2d^2\underline{l})=0\,\forall\, i\geq 1\iff h^i_{b_X^*\mathcal{E}\otimes N^{\otimes ab}}((r/s+c/d-i)b^2\underline{l})=0\,\forall\, i\geq 1.
\end{equation}
But the left hand side of \eqref{q0} is equivalent to the following condition
\begin{equation}\label{q0'}
    h^i_{\textrm{gen}}((bds)_X^*\mathcal{E}\otimes N^{\otimes abd^2s^2}\otimes L^{\otimes b^2cds^2}\otimes L^{\otimes (r-is)b^2d^2s})=0\,\forall\, i\geq 1.
\end{equation}
And the right hand side of \eqref{q0} is equivalent to $h^i_{\textrm{gen}}((bds)_X^*\mathcal{E}\otimes N^{\otimes abd^2s^2}\otimes L^{\otimes (rd+cs-ids)b^2ds})=0\,\forall\, i\geq 1$ which under simplification boils down to \eqref{q0'}. The proof is now complete.\QEDB

\begin{proposition}\label{cont2}
Let $X$ be an abelian variety of dimension $g$ and let $\mathcal{E}$ be a vector bundle on $X$. Let $\underline{l},\underline{l'},\underline{n}\in N^1(X)$ with $\underline{l},\underline{l'}$ ample. Further, let $x=a/b,y=c/d\in\mathbb{Q}$ with $a,b,c,d\in\mathbb{Z}$, $b,c,d>0$. Then $$\mathbb{Q}\textrm{-reg}_{\underline{l}}(\mathcal{E}\langle x\underline{n}+y\underline{l'}\rangle)\leq\mathbb{Q}\textrm{-reg}_{\underline{l}}(\mathcal{E}\langle x\underline{n}\rangle).$$
\end{proposition}

\noindent\textit{Proof.} Let $\beta=\beta_1/\beta_2\in\mathbb{Q}$ with $\beta_1,\beta_2\in\mathbb{Z}$, $\beta_2>0$. It is enough to show that if $\mathbb{Q}\textrm{-reg}_{\underline{l}}(\mathcal{E}\langle x\underline{n}\rangle)<\beta$ then $\mathbb{Q}\textrm{-reg}_{\underline{l}}(\mathcal{E}\langle x\underline{n}+y\underline{l'}\rangle)\leq \beta$. As before, $\mathbb{Q}\textrm{-reg}_{\underline{l}}(\mathcal{E}\langle x\underline{n}\rangle)<\beta\implies h^i_{b_X^*\mathcal{E}\otimes N^{\otimes ab}}((\beta-i)b^2\underline{l})=0\,\forall\,i\geq 1$ which is equivalent to $(b\beta_2)_X^*\mathcal{E}\otimes N^{\otimes ab\beta_2^2}\otimes L^{\otimes\beta_1\beta_2b^2}$ being continuously $0$--regular for $(X,L^{\otimes (b\beta_2)^2})$. This in turn is equivalent to $\mathcal{G}:=(bd\beta_2)_X^*\mathcal{E}\otimes N^{\otimes abd^2\beta_2^2}\otimes L^{\otimes \beta_1\beta_2b^2d^2}$ being continuously $0$--regular for $(X,L^{\otimes (bd\beta_2)^2})$. 

We aim to show that $\mathbb{Q}\textrm{-reg}_{\underline{l}}(\mathcal{E}\langle x\underline{n}+y\underline{l'}\rangle)\leq\beta$ which is equivalent to showing the vanishing condition $h^i_{(bd)_X^*\mathcal{E}\otimes N^{\otimes abd^2}\otimes L'^{\otimes b^2cd}}((\beta-i)b^2d^2\underline{l})=0\,\forall\,i\geq 1$ which in turn is equivalent to showing that $\mathcal{G}\otimes L'^{\otimes b^2cd\beta_2^2}$ is continuously $0$--regular for $(X,L^{\otimes (bd\beta_2)^2})$. But this follows from Corollary \ref{earlier}.\QEDB

\vspace{5pt}

\noindent\textit{Proof of Theorem \hyperref[qregthm]{B}.} We adapt an argument of Ito (proof of \cite[Proposition 2.9]{Ito}). First we show that $\mathbb{Q}\textrm{-reg}_{\underline{l}}(\mathcal{E}\langle - \rangle): N^1(X)_{\mathbb{Q}}\to\mathbb{R}$ is continuous. Let $\left\{\xi_i\right\}_{i\in\mathbb{N}}$ be a sequence in $N^1(X)_{\mathbb{Q}}$ converging to $\xi\in N^1(X)_{\mathbb{Q}}$. For any rational number $\delta>0$, there exists $N_0\in\mathbb{N}$ such that for all $i\geq N_0$, $\xi_{i}-\xi+\delta\underline{l}$ and $\xi+\delta\underline{l}-\xi_i$ are both ample. Thus, by Proposition \ref{cont2}, we deduce that for all $i\geq N_0$, we have the following inequality
\begin{equation*}
    \mathbb{Q}\textrm{-reg}_{\underline{l}}(\mathcal{E}\langle \xi+\delta\underline{l}\rangle)\leq \mathbb{Q}\textrm{-reg}_{\underline{l}}(\mathcal{E}\langle \xi_i\rangle)\leq \mathbb{Q}\textrm{-reg}_{\underline{l}}(\mathcal{E}\langle \xi-\delta\underline{l}\rangle).
\end{equation*}
Now by Lemma \ref{cont1}, we deduce that for all $i\geq N_0$, $|\mathbb{Q}\textrm{-reg}_{\underline{l}}(\mathcal{E}\langle \xi_i\rangle)-\mathbb{Q}\textrm{-reg}_{\underline{l}}(\mathcal{E}\langle \xi\rangle)|\leq \delta$ and that proves the claim. To finish the proof, we need to show that given a sequence $\left\{\xi_i\right\}_{i\in\mathbb{N}}$ in $N^1(X)_{\mathbb{Q}}$ converging to $\xi\in N^1(X)_{\mathbb{R}}$, $\left\{\mathbb{Q}\textrm{-reg}_{\underline{l}}(\mathcal{E}\langle \xi_i\rangle)\right\}_{i\in\mathbb{N}}$ converges to a real number. As before, for any rational $\delta>0$, there exists $N_0\in\mathbb{N}$ such that for $j,k\geq N_0$, $\xi_j-\xi+\delta\underline{l}$ and $\xi+\delta\underline{l}-\xi_k$ are both ample. Using Proposition \ref{cont2}, we deduce that $\mathbb{Q}\textrm{-reg}_{\underline{l}}(\mathcal{E}\langle \xi_j+\delta\underline{l}\rangle)\leq \mathbb{Q}\textrm{-reg}_{\underline{l}}(\mathcal{E}\langle \xi_k-\delta\underline{l}\rangle)$ which by Lemma \ref{cont1} implies that for all $j,k\geq N_0$, we have the inequality $\mathbb{Q}\textrm{-reg}_{\underline{l}}(\mathcal{E}\langle \xi_j\rangle)-\delta\leq \mathbb{Q}\textrm{-reg}_{\underline{l}}(\mathcal{E}\langle \xi_k\rangle)+\delta$.  Thus we obtain
\begin{equation*}
    \uplim_{j\to\infty} \mathbb{Q}\textrm{-reg}_{\underline{l}}(\mathcal{E}\langle \xi_j\rangle)-\delta\leq \lowlim_{k\to\infty} \mathbb{Q}\textrm{-reg}_{\underline{l}}(\mathcal{E}\langle \xi_k\rangle)+\delta.
\end{equation*}
Now, by letting $\delta\to 0$, we see that $\left\{\mathbb{Q}\textrm{-reg}_{\underline{l}}(\mathcal{E}\langle \xi_i\rangle)\right\}_{i\in\mathbb{N}}$ converges to a real number.\QEDB

\section{Applications for sheaves on abelian varieties}\label{seccor} 

\subsection{Three immediate consequences} We list three consequences of our generic vanishing theorems.

\begin{corollary}\label{v0}
Let $(X,H)$ be a polarized smooth projective variety. Assume that there exists a globally generated line bundle $H_1$ on $X$ and an ample line bundle $H_2$ on $\textrm{Alb}(X)$ such that $H=H_1+\textrm{alb}_X^*H_2$, where $\textrm{alb}_X:X\to\textrm{Alb}(X)$ is the Albanese map. Let $\mathcal{F}$ be a torsion--free coherent sheaf on $X$ that is continuously $1$--regular for $(X,H)$. Then $\chi(\mathcal{F})\geq 0$ with equality if and only if either $V^0(\mathcal{F})=\emptyset$ or any component $V^0(\mathcal{F})$ is of codimension $1$.
\end{corollary}

\noindent\textit{Proof.} We know that $\mathcal{F}$ is GV from Theorem \hyperref[newgv]{A}, and moreover $V^i(\mathcal{F})=\emptyset$ for $i\geq 2$. It follows from the generic vanishing theory that $\chi(\mathcal{F})\geq 0$ with equality if and only if $\textrm{codim}(V^0(\mathcal{F}))\geq 1$. Now, by \cite[Proposition 3.15]{PaPo4}, and \cite[Lemma 1.8]{Par'}, we know that a codimension $q$ component of $V^0(\mathcal{F})$ is also a component of $V^q(\mathcal{F})$ whence the conclusion follows.\QEDB

\begin{remark}
For a coherent sheaf $\mathcal{F}$ on a smooth projective variety $X$ of dimension $d$, define $R\Delta\mathcal{F}:={\bf R}\mathcal{H}\textrm{om}(\mathcal{F},K_X)$. If $\mathcal{F}$ is GV, then ${\bf R}\Phi_{P}(R\Delta\mathcal{F})$ is a sheaf concentrated in degree $d$, i.e., ${\bf R}\Phi_P(R\Delta\mathcal{F})=R^d\Phi_P(R\Delta\mathcal{F})[-d]$. The generic vanishing theory tells us that for a GV sheaf $\mathcal{F}$, $\chi(\mathcal{F})$ is the rank of $\widehat{R\Delta\mathcal{F}}$ and the support of $\widehat{R\Delta\mathcal{F}}$ is $-V^0(\mathcal{F})$ where $\widehat{R\Delta\mathcal{F}}=R^d\Phi_P(R\Delta\mathcal{F})$. Thus, if $X, H, \mathcal{F}$ are as in Corollary \ref{v0}, then $\chi(\mathcal{F})\geq 0$ with equality if and only if any (non--empty) component of the support of $\widehat{R\Delta\mathcal{F}}$ is of codimension one.
\end{remark}

It is worth mentioning that for a GV sheaf $\mathcal{F}$ on an abelian variety, $V^0(\mathcal{F})\neq \emptyset$ by \cite[Lemma 1.12]{Par'}.

\begin{remark}
The conclusion of Corollary \ref{v0} also applies in the set--up of Corollaries \ref{newgvcor}, \ref{var1} and Theorem \ref{var2}, all with $k=1$.
\end{remark}

\begin{corollary}\label{gvnef}
Let $(X,H)$ be a polarized abelian variety and let $\mathcal{F}$ be torsion--free coherent sheaf on $X$. Assume $\mathcal{F}$ is continuously $k$--regular for $(X,H)$ for some $k\in\mathbb{N}$. Then the following statements hold.
\begin{itemize}
    \item[(1)] If $1\leq k\leq \dim X$, then $\mathcal{F}$ is a GV$_{-(k-1)}$ sheaf. In particular, if $k=1$, then $\mathcal{F}$ is nef.
    \item[(2)] If $k=0$, then $\mathcal{F}$ is an IT$_0$ sheaf, in particular $\mathcal{F}$ is ample.
\end{itemize}
\end{corollary}

\noindent\textit{Proof.} Follows immediately from Theorem \hyperref[newgv]{A} combined with the facts that GV sheaves on abelian varieties are nef (\cite[Theorem 4.1]{PaPo3}), and IT$_0$ sheaves are ample (\cite[Corollary 3.2]{Deb}).\QEDB

\begin{remark}
Let $(X,H)$ be a polarized abelian variety and let $\mathcal{E}$ be a vector bundle of rank $r$ that is continuously $0$--regular for $(X,H)$. Then for any subvariety $Z\subseteq X$ of dimension $k$, we have $$c_1(\mathcal{E})^k\cdot Z=(c_1(\mathcal{E}(-H))+rH)^k\cdot Z\geq r^kH^kZ$$ as $\mathcal{E}(-H)$ is nef by the above corollary. In particular, $c_1(\mathcal{E})^{n}\geq r^nH^n$ where $n:=\dim X$, and for any $x\in X$, we have the Seshadri constant $\epsilon(\textrm{det}(\mathcal{E}),x)\geq r\epsilon(H,x)$. This refines the bound of \cite[Theorem 7.2]{Lop1} for abelian varieties.
\end{remark}

\begin{corollary}\label{regin}
Let $(X,H)$ be a polarized abelian variety. Let $\mathcal{F}_1$ and $\mathcal{F}_2$ be torsion--free coherent sheaves on $X$ with one of them locally free. Assume $\mathcal{F}_1$ is continuously $k_1$--regular and $\mathcal{F}_2$ is continuously $k_2$--regular for $(X,H)$. Then $\mathcal{F}_1\otimes\mathcal{F}_2$ is $C_{0,k_1+k_2}$ for $(X,H)$.
\end{corollary}

\noindent\textit{Proof.} Without loss of generality, we may assume that $k_1=k_2=0$. For any integer $1\leq i\leq \dim X$, we have  $H^i(\mathcal{F}_1\otimes\mathcal{F}_2(-iH))=H^i((\mathcal{F}_1(-iH))\otimes \mathcal{F}_2)$. By Theorem \hyperref[newgv]{A}, we know that $\mathcal{F}_1(-iH)$ is GV$_{-(i-1)}$ and $\mathcal{F}_2(-H)$ is GV. By the preservation of vanishing, $\mathcal{F}_2$ is IT$_0$, whence by Proposition \ref{preservation}, we obtain the required vanishing $H^i(\mathcal{F}_1\otimes\mathcal{F}_2(-iH))=0$.\QEDB

\begin{remark}
We remark that the proof of the above corollary also shows the following: let $(X,H)$ be a polarized abelian variety and let $\mathcal{F}_1$, $\mathcal{F}_2$ be torsion--free coherent sheaves on $X$ with one of them is locally free. If $\mathcal{F}_1$ is continuously $0$--regular for $(X,H)$ and $\mathcal{F}_2$ is IT$_0$, then $\mathcal{F}_1\otimes\mathcal{F}_2$ is $0$--regular for $(X,H)$. This can also be thought of as a generalization of a result of Murty and Sastry (\cite[Proposition 5.4.1]{MS}). 
\end{remark}

\subsection{Syzygies of tautological bundles of zero--regular bundles}\label{5} We start with some background on syzygies.

\subsubsection{} Let $X$ be a smooth projective variety and let $L$ be a very ample line bundle on $X$. Consider the embedding $X\subseteq \mathbb{P}^r$ given by the complete linear series $|L|$ where $r=h^0(L)-1$. One has a minimal graded free resolution of $R(X,L):=\bigoplus\limits_{q\geq 0} H^0(qL)$ as an $S:=\textrm{Sym}^{\bullet}(H^0(L))$ module as follows:
\begin{equation*}
    0\to E_{r+1}=\bigoplus\limits_{j}S(-a_{r+1,j})\to E_r=\bigoplus\limits_{j}S(-a_{r,j})\to\cdots\to E_1=\bigoplus\limits_{j}S(-a_{1,j})\to E_0=\bigoplus\limits_{j}S(-a_{0,j})\to R(X,L)\to 0.
\end{equation*}
For a reference of the following definition, see for example \cite[Definition 1.8.50]{Laz}.

\begin{definition} {\bf (Projective normality and $N_p$ property)} Suppose we are in the situation as above.
\begin{itemize}
    \item[--] The embedding given by the complete linear series $|L|$ is called {\it projectively normal} if $E_0=S$.
    \item[--] We say that {\it $L$ satisfies $N_p$ property} if the embedding by $|L|$ is projectively normal, and $a_{ij}=i+1$ for all $1\leq i\leq p$.
\end{itemize}
\end{definition}

In practice, to calculate the syzygies of a projective variety one needs to calculate cohomology groups involving the syzygy bundles that we will define next (see for example \cite[Section 3]{Park}).

\begin{definition}
{\bf (Syzygy bundle)} Let $X$ be a smooth projective variety and let $\mathcal{E}$ be a globally generated vector bundle on $X$. The {\it syzygy bundle} $M_{\mathcal{E}}$ is the kernel of the map $H^0(\mathcal{E})\otimes \mathcal{O}_X \to \mathcal{E}$ i.e. we have the exact sequence
\begin{equation}\label{kernel}
    0 \to M_{\mathcal{E}} \to  H^0(\mathcal{E})\otimes \mathcal{O}_X \to \mathcal{E} \to 0.
\end{equation}
\end{definition}

The following proposition of Park will be used in the proof of Theorem \ref{thm_appendix} and we include it here. 

\begin{proposition}\label{parks}
(\cite[Proposition 3.2]{Park}) Let $X$ be a smooth projective variety and let $\mathcal{E}$ be an ample and globally generated bundle on $X$. Then $\mathcal{O}_{\mathbb{P}(\mathcal{E})}(1)$ on $\mathbb{P}(\mathcal{E})$ satisfies $N_p$ property if $H^k\left(\bigwedge\limits^{i}M_{\mathcal{E}}\otimes\mathcal{E}^{\otimes j}\right)=0$ for all $0\leq i\leq p+1$, $j,k\geq 1$.
\end{proposition} 

We now prove Corollary \hyperref[npmain]{C} that is obtained by an immediate application of our result combined with the theorem of Ito (Theorem \ref{thm_appendix}).

\vspace{5pt}

\noindent\textit{Proof of Corollary \hyperref[npmain]{C}.} Notice that for any $x\in\mathbb{Q}$ with $x<1$, we have $\mathbb{Q}\textrm{-reg}_{\underline{h}}(\mathcal{E}\langle -x\underline{h}\rangle)<1$ by Lemma \ref{cont1}, whence $\mathcal{E}\langle-x\underline{h}\rangle$ is IT$_0$ by Corollary \ref{qtwist} . Since $\beta(\underline{h})<\frac{1}{p+2}$, we have $(p+2)\beta(\underline{h})<1$. Thus, for $0<\epsilon\ll 1$, we have $x=1-\epsilon\geq (p+2)\beta(\underline{h})$ whence the assertion follows immediately by Theorem \ref{thm_appendix}.\QEDB

\appendix

\section{On Syzygies of projective bundles on abelian varieties \texorpdfstring{\\}{TEXT} by Atsushi Ito}\label{appendix}

\vspace{3mm}

In this appendix, we follows the notation in the previous sections.
In particular, we work over the field of complex numbers $\mathbb{C}$.

Let $(X,L)$ be a polarized abelian variety and  $\underline{l} \in N^1(X) $ be the class of $L$.
In \cite{JP}, Z.\ Jiang and G.\ Pareschi define the \emph{basepoint-freeness threshold} $\beta(\underline{l}) \in (0,1]$.
This invariant is quite useful to study syzygies of polarized abelian varieties
since \cite{JP} and F.~Caucci \cite{Cau} show that $L$ satisfies $N_p$ property if $\beta(\underline{l}) < 1/(p+2)$.

The purpose of this short note is to prove the following theorem:

\begin{theorem}\label{thm_appendix}
Let $(X,L)$ be a polarized abelian variety and $\mathscr{E}$ be a vector bundle on $ X$.
Let $p \geq 0$ be an integer.
Assume that there exists a rational number $x \geq (p+2) \beta(\underline{l})$ such that
$\mathscr{E} \langle  -x  \underline{l} \rangle$ is M-regular.
Then $\mathscr{O}_{\mathbb{P}(\mathscr{E})}(1)$ satisfies $N_p$ property.
\end{theorem}

\begin{remark}
In the case $\mathscr{E}=L$,
$L \langle  -x  \underline{l} \rangle $ is M-regular if and only if $ x <1$ (cf.\ \cite[Example 2.1]{Ito22}).
Hence the existence of a rational number $x$ in Theorem \ref{thm_appendix} is equivalent to $\beta (\underline{l})< 1/(p+2)$ in this case.
\end{remark}

To prove this theorem, we use the following lemma.

\begin{lemma}\label{lem_appendix}
Let $(X,L)$ be a polarized abelian variety and $\mathscr{E}$ be a vector bundle on $ X$.
Assume that there exists a rational number $x \geq \beta(\underline{l})$ such that
$\mathscr{E} \langle  -x  \underline{l} \rangle$ is M-regular. Then
\begin{enumerate}
\item $ \mathscr{E}$ is IT$_0$ and globally generated.
\item Let $M_{\mathscr{E}}$ be the  syzygy bundle of $\mathscr{E}$
defined by \eqref{kernel}.
For a rational number $y >0 $, $M_{\mathscr{E}} \langle y \underline{l} \rangle$ is IT$_0$
if $\tfrac1x + \tfrac1y  \leq \tfrac1{\beta(\underline{l}) }$.
\end{enumerate}
\end{lemma}

\noindent\textit{Proof.}
(1) Since $\mathscr{E} \langle  -x  \underline{l} \rangle$ is M-regular and  $x \geq \beta(\underline{l}) >0$,
 the bundle $ \mathscr{E}$ is IT$_0$.
The global generation of $\mathscr{E}$ follows from \cite[Theorem 1.2 (1)]{Ito22}.\\
(2) Write $y=a/b $ with integers $ a, b >0$.
Then 
\begin{align*}
\text{$M_{\mathscr{E}} \langle y \underline{l} \rangle =M_{\mathscr{E}} \langle \tfrac{a}{b} \underline{l} \rangle$ is IT$_0$} \quad  
&\iff \quad  \text{$b_X^* M_{\mathscr{E}}  \otimes L^{\otimes ab}$ is IT$_0$} \\
&\iff \quad \text{${b_X}_{*} (b_X^* M_{\mathscr{E}}  \otimes L^{\otimes ab}) = M_{\mathscr{E}}  \otimes {b_X}_{*} (L^{\otimes ab})$ is IT$_0$},
\end{align*}
where the first equivalence follows from the definition
and 
the second one holds since the property IT$_0$ is preserved by the pushforward by isogenies (cf.\ \cite[(2.4)]{Ito22}).

Consider the exact sequence
\[
0 \rightarrow M_{\mathscr{E}} \otimes {b_X}_{*} (L^{\otimes ab}) \rightarrow H^0(\mathscr{E}) \otimes {b_X}_{*} (L^{\otimes ab}) \rightarrow \mathscr{E} \otimes {b_X}_{*} (L^{\otimes ab}) \rightarrow 0
\]
obtained by tensoring ${b_X}_{*} (L^{\otimes ab})$ with  \eqref{kernel}.
Since $L^{\otimes ab}$ is IT$_0$,
so is the pushforward ${b_X}_{*} (L^{\otimes ab})$.
Since $\mathscr{E}$ is  IT$_0$ as well by (1),
both $H^0(\mathscr{E}) \otimes {b_X}_{*} (L^{\otimes ab}) $ and $\mathscr{E} \otimes {b_X}_{*} (L^{\otimes ab})$ are  IT$_0$.
Thus $H^i(M_{\mathscr{E}} \otimes {b_X}_{*} (L^{\otimes ab})  \otimes \zeta)  =0$ for any $i \geq 2$ and $\zeta \in \textrm{Pic}^0 (X)$,
and 
hence $M_{\mathscr{E}} \otimes {b_X}_{*} (L^{\otimes ab}) $ is  IT$_0$ if and only if 
\[
H^1(M_{\mathscr{E}} \otimes {b_X}_{*} (L^{\otimes ab})  \otimes \zeta) =0
\]
for any $\zeta \in \textrm{Pic}^0 (X)$,
which is equivalent to the surjectivity of
\begin{align}\label{eq_appendix2}
H^0(\mathscr{E}) \otimes H^0({b_X}_{*} (L^{\otimes ab})  \otimes \zeta) \rightarrow H^0(\mathscr{E} \otimes {b_X}_{*} (L^{\otimes ab}) \otimes \zeta).
\end{align}

In conclusion, (2) follows from the surjectivity of \eqref{eq_appendix2} for any $\zeta \in \textrm{Pic}^0 (X)$.
Set $\mathscr{F} = {b_X}_{*} (L^{\otimes ab})  \otimes \zeta$.
Recall that $y = a/b$.
Since
\[
\frac1x+\frac{b}{a} = \frac1x+\frac1y \leq 
\frac1{\beta(\underline{l}) }
\]
by assumption,
it suffices to show that
\[
\varphi_{\underline{l}}^* \Phi ( \mathscr{E}) \otimes \varphi_{\underline{l}}^* \Phi ( (-1)_X^* \mathscr{F}) \langle (\tfrac1x+\tfrac{b}{a}) \underline{l} \rangle
\]
is M-regular by \cite[Proposition 4.4]{Ito22},
where $\varphi_{\underline{l}}  : X \rightarrow \hat{X}=\textrm{Pic}^0(X)$ is the isogeny induced by the polarization $\underline{l}$,
$ \Phi =\mathbf{R} \hat{\mathscr{S}} : \mathbf{D}(X) \rightarrow \mathbf{D}(\hat{X})$ is the Fourier-Mukai functor associated to the Poincar\'e line bundle on $X \times \hat{X}$,
 and $(-1)_X$ is the multiplication map by $(-1)$ on $X$.
We note that $\Phi(\mathscr{E}), \Phi((-1)_X^* \mathscr{F})$ are locally free sheaves since $\mathscr{E}, (-1)_X^* \mathscr{F}$ are IT$_0$.

The rest is to check the M-regularity of $\varphi_{\underline{l}}^* \Phi ( \mathscr{E}) \otimes \varphi_{\underline{l}}^* \Phi ( (-1)_X^* \mathscr{F}) \langle ( \tfrac1x+\tfrac{b}{a} ) \underline{l}\rangle $,
which can be shown by modifying the argument in the proof of \cite[Proposition 4.4 (1)]{Ito22} as follows:
The pullback of a $\mathbb{Q}$-twisted sheaf $\mathscr{G} \langle t \underline{l} \rangle $ by the multiplication map $n_X $ on $X$ for $n \geq 1$ is defined as
\[
 n_X^* (\mathscr{G} \langle t \underline{l} \rangle) := n_X^* \mathscr{G} \langle t  \,  n_X^*\underline{l} \rangle = n_X^* \mathscr{G} \langle n^2 t \underline{l} \rangle  .
\]
Since M-regularity is preserved by such pullbacks (cf.\ \cite[(2.2)]{Ito22}),
it suffices to show the M-regularity of
\begin{align}\label{eq_appendix3}
a_X^* \left( \varphi_{\underline{l}}^* \Phi ( \mathscr{E}) \otimes \varphi_{\underline{l}}^* \Phi ( (-1)_X^* \mathscr{F}) \langle (\tfrac1x+\tfrac{b}{a})\underline{l} \rangle  \right)
= a_X^* \left( \varphi_{\underline{l}}^* \Phi ( \mathscr{E}) \otimes \varphi_{\underline{l}}^* \Phi ((-1)_X^*  \mathscr{F}) \right)  \langle (\tfrac{a^2}x+ab) \underline{l} \rangle  .
\end{align}
Recall that $\mathscr{F} = {b_X}_{*} (L^{\otimes ab})  \otimes \zeta = {b_X}_{*} (L^{\otimes ab} \otimes b_X^* \zeta)   $.
We take $\zeta' \in \textrm{Pic}^0  (X) $ such that $\zeta'^{\otimes ab} = b_X^* \zeta $
and set $L'= (-1)_X^* (L \otimes \zeta')$.
Then we have
\[
(-1)_X^*\mathscr{F} =(-1)_X^*  {b_X}_* ( (L \otimes \zeta')^{\otimes ab} )= {b_X}_* ((-1)_X^*  (L \otimes \zeta')^{\otimes ab} )= {b_X}_* ( L'^{\otimes ab})
\]
and hence
\begin{align*}
 a_X^*  \varphi_{\underline{l}}^* \Phi ((-1)_X^* \mathscr{F}) =  a_X^*  \varphi_{\underline{l}}^* \Phi ( {b_X}_* ( L'^{\otimes ab}) )
 =a_X^*  \varphi_{\underline{l}}^* b_{\hat{X}}^*  \Phi ( L'^{\otimes ab} )  =  \varphi_{ab\underline{ l}}^* \Phi ( L'^{\otimes ab} )
 =H^0(L'^{\otimes ab} ) \otimes L'^{\otimes -ab}.
\end{align*}
In fact, the second equality follows from \cite[(3.4)]{Muk81}.
The third one follows from $ b_{\hat{X}} \circ \varphi_{\underline{l}} \circ a_X = \varphi_{ ab \underline{l}} $.
The last one follows from \cite[Proposition 3.11 (1)]{Muk81}
since the numerical class of  $ L'^{\otimes ab}$ is $ ab \underline{l}$ by $L' = (-1)_X^* (L \otimes \zeta')\equiv L$.
Hence it holds that
\begin{align*}
a_X^* \left( \varphi_{\underline{l}}^* \Phi ( \mathscr{E}) \otimes \varphi_{\underline{l}}^* \Phi ( (-1)_X^*\mathscr{F}) \right)  \langle (\tfrac{a^2}x+ab) \underline{l} \rangle  
&= a_X^*  \varphi_{\underline{l}}^* \Phi ( \mathscr{E}) \otimes H^0(L'^{\otimes ab} ) \otimes L'^{\otimes -ab} \langle (\tfrac{a^2}x+ab)\underline{l} \rangle \\
&= a_X^*  \varphi_{\underline{l}}^* \Phi ( \mathscr{E}) \otimes H^0(L'^{\otimes ab} )  \langle \tfrac{a^2}x \underline{l} \rangle\\
&= H^0(L'^{\otimes ab})  \otimes a_X^*  \left( \varphi_{\underline{l}}^* \Phi ( \mathscr{E}) \langle \tfrac1x \underline{l} \rangle \right),
\end{align*}
where the second equality follows from $L' \equiv L$.
Since $\mathscr{E} \langle  -x  \underline{l} \rangle$ is M-regular,
$\varphi_{\underline{l}}^* \Phi ( \mathscr{E}) \langle \tfrac1x \underline{l} \rangle$ is M-regular as well by \cite[Proposition 4.1]{Ito22}.
Thus so is the pullback $a_X^*  \left( \varphi_{\underline{l}}^* \Phi ( \mathscr{E}) \langle \tfrac1x \underline{l} \rangle \right)$
and hence we obtain the M-regularity of \eqref{eq_appendix3},
which implies (2).\QEDB

\vspace{5pt}

\noindent\textit{Proof of Theorem \ref{thm_appendix}.}
Since  $x \geq (p+2) \beta(\underline{l}) \geq  \beta(\underline{l}) $, the bundle $\mathscr{E}$ is globally generated by Lemma \ref{lem_appendix} (1).

To prove $N_p$ property for $\mathscr{O}_{\mathbb{P}(\mathscr{E})}(1)$,
it suffices to show that
$H^k (M_{\mathscr{E}}^{\otimes i} \otimes \mathscr{E}^{\otimes j}) =0$ for all $ 0 \leq i \leq p+1, j,k \geq 1$ by Proposition \ref{parks}
since we work over $\mathbb{C}$, whose characteristic is zero.
Hence this theorem holds if $M_{\mathscr{E}}^{\otimes i} \otimes \mathscr{E}^{\otimes j}$ is IT$_0$ for all $ 0 \leq i \leq p+1, j \geq 1$.

If $i=0$, $M_{\mathscr{E}}^{\otimes i} \otimes \mathscr{E}^{\otimes j} = \mathscr{E}^{\otimes j}$ is IT$_0$  for $j \geq 1$
since so is $\mathscr{E}$ by Lemma \ref{lem_appendix} (1).
If $i \geq 1$,
 $M_{\mathscr{E}}^{\otimes i} \otimes \mathscr{E}^{\otimes j} $ is written as
\begin{align}\label{eq_appendix_M_E}
M_{\mathscr{E}}^{\otimes i} \otimes \mathscr{E}^{\otimes j} = 
\left(M_{\mathscr{E}} \langle  \tfrac{x}{i } \underline{l} \rangle \right)^{\otimes i} \otimes \mathscr{E} \langle - x \underline{l} \rangle \otimes \mathscr{E}^{\otimes j-1}
\end{align}
as a $\mathbb{Q}$-twisted sheaf.
For  $1 \leq i \leq p+1$, $M_{\mathscr{E}} \langle  \tfrac{x}{i } \underline{l} \rangle$ is IT$_0$ by Lemma \ref{lem_appendix} (2) 
since it holds that
\[
\frac1x + \frac{i}x \leq  \frac{p+2}{x} \leq \frac1{\beta(\underline{l}) }
\]
by assumption.
Furthermore, $ \mathscr{E} \langle - x \underline{l} \rangle$ is M-regular by assumption and $\mathscr{E}$ is IT$_0$.
Thus their tensor product \eqref{eq_appendix_M_E} is IT$_0$ by \cite[Proposition 3.4]{Cau}.\QEDB

\subsection*{Acknowledgments}
A.I.\ was supported by JSPS KAKENHI Grant Numbers 17K14162, 21K03201.

\bibliographystyle{plain}

\end{document}